\documentclass[12pt]{article}
\usepackage[utf8]{inputenc}

\usepackage{style}

\title{Decoupled Solution for Composite Sparse-plus-Smooth Inverse Problems\footnote{The associated code is freely accessible on the dedicated open-access repository: \url{https://github.com/AdriaJ/decouple-composite}}}
\author{Adrian Jarret$^1$\thanks{ Corresponding author,\url{adrian.jarret@epfl.ch}.} \and Julien Fageot$^1$}

\date{
	$^1${\small Laboratoire des Communications Audiovisuelles,\\ \'Ecole Polytechnique Fédérale de Lausanne,\\ 1015 Lausanne, Switzerland}%
}

\begin{document}

\maketitle

\begin{abstract}
    We consider composite linear inverse problems where the signal to recover is modeled as a superposition of two functions. 
    Relying on a variational framework, we formulate an optimization problem over the pairs of components using two regularization terms, each of them acting on a different part of the solution.
    One component is regularized with a quadratic norm over a Hilbert space, which promotes a smooth solution. The second component is less constrained and is modeled as a Banach space element, allowing for the promotion of sparsity.

    We show how this composite optimization problem can be reduced to an optimization problem over the Banach space component only up to a linear problem.
    This reveals a decoupling between the two components, leading to a new composite representer theorem.
    It naturally induces a decoupled numerical procedure to solve the composite optimization problem.

    We exemplify our main result with a composite deconvolution problem of Dirac recovery over a smooth background. In this setting, we illustrate the relevance of a composite model and show a significant temporal gain on signal reconstruction, which results from our decoupled algorithmic approach.
\end{abstract}


\section{Introduction}


    \subsection{Composite linear inverse problems}
    
\begin{figure}[t]
        \centering
        \includegraphics[width=0.8\linewidth]{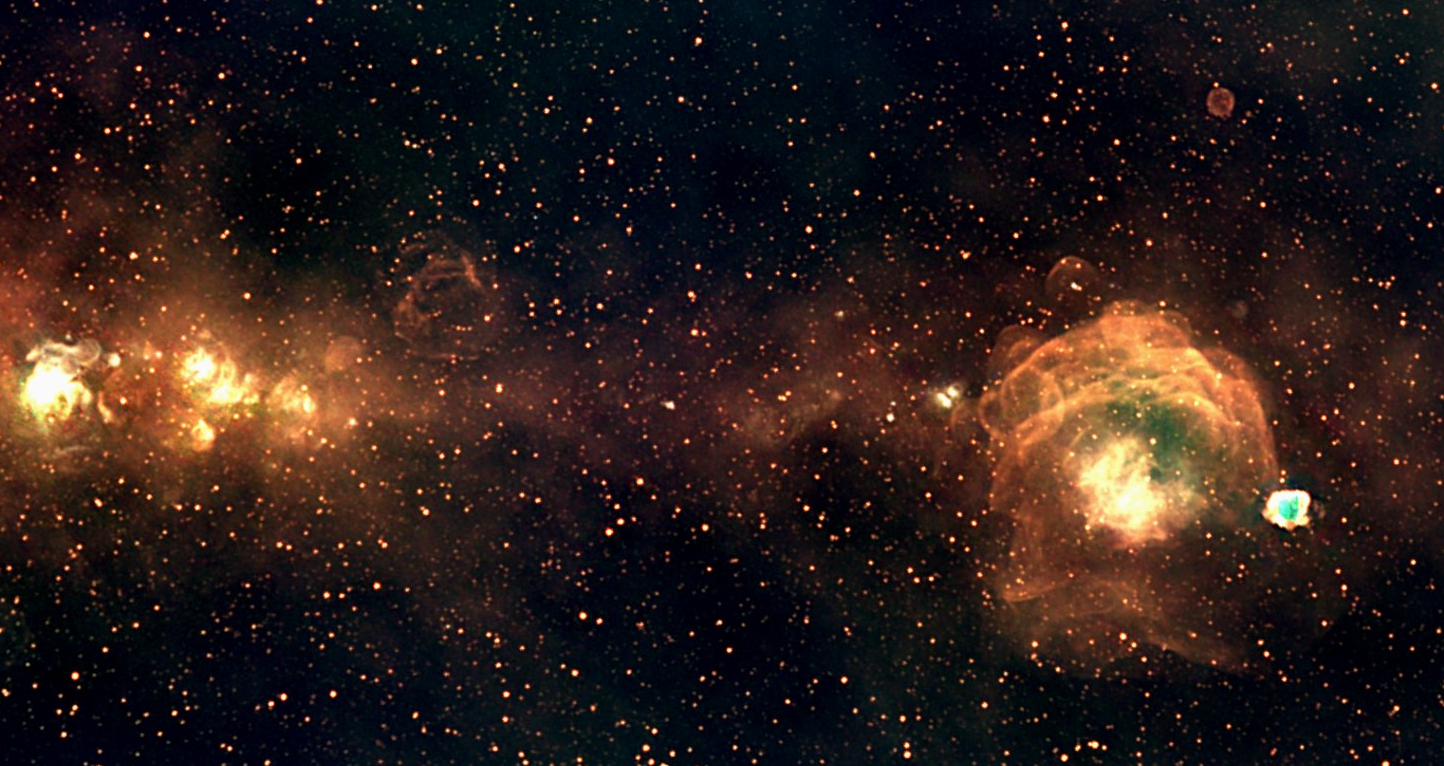}
        \caption{Observation of the radio sky at coordinates from the GLEAM survey accessible at \href{https://gleamoscope.icrar.org/gleamoscope/trunk/src/}{gleamoscope.icrar.org}, J2000 coordinates (9h37min15.21s, 50\textdegree25'03.1'').}
        \label{fig:gleam-sky}
    \end{figure}

    We study linear inverse problems where the goal is to recover a signal $s^\dagger$ from a finite number of possibly noisy linear measurements $\bm{y} \approx \bm{\Phi}(s^\dagger)$. 
    The inverse problem is called composite if the signal to reconstruct is modeled as a sum of several subcomponents. We consider the case where $s^\dagger$ is the sum of two terms
    \begin{equation}
        s^\dagger = s_1^\dagger + s_2^\dagger
        \label{eq:summodel}        
    \end{equation}
    with different characteristics.
    Such a framework can be used for instance to model signals which contain mixed-types information, such as background and foreground in imaging. A concrete example of composite-type signal is provided in figure~\ref{fig:gleam-sky} with an image of the sky from radio observations, containing both precisely-localized point sources (stars) and diffuse emissions (nebulae) \cite{hurley-walker2017}.
    Our goal is to approximately recover not only $s^\dagger$ but also the subcomponents $s_1^\dagger$ and $s_2^\dagger$ from the observations $\bm{y}$.

    In a theoretical model, the addition between $s_1^\dagger$ and $s_2^\dagger$ may not be properly defined because the signals may belong to different spaces. For this reason, in what follows we consider the more general inverse problem over the pair $(s_1^\dagger, s_2^\dagger)$ with two possibly different measurement operators
    \begin{equation}
        \label{eq:generic-comp-ip}
        \bm{y} \approx \bm{\Phi}_1(s_1^\dagger) + \bm{\Phi}_2(s_2^\dagger)
    \end{equation}

    The inverse problem \eqref{eq:generic-comp-ip} is generally ill-posed. First, the data are typically insufficient to reconstruct the total information of the original signal: many signals can explain the observation vector $\bm{y}$. Second, the observation vector $\bm{y}$ is often a noisy version of the ideal measurement vector $\bm{\Phi}_1(s_1^\dagger) + \bm{\Phi}_2(s_2^\dagger)$. Third, the search for the subcomponents comes with an additional difficulty since it requires to separate the sources explaining adequately the observations.

    Various classical signal recovery applications are based on a composite model for linear inverse problems.
    For instance, structure-texture decomposition in images corresponds to a simple inverse problem (no forward operator) for which numerous composite models exist \cite{aujol2006std,guennec2024std}.
    We can also refer to background-foreground separation in images, which usually relies on composite models. They include \emph{sparse-plus-smooth} models, for instance used for deconvolution imaging in microscopy \cite{debarnot2020learning,stergiopoulou2022colorme}, or \emph{low-rank-plus-sparse} models for detecting moving objects from a background in video analysis \cite{bouwmans2017decomposition}, for which theoretic reconstruction guarantees exist \cite{tanner2023composite}.
    Composite modeling has also proven relevant for other types of signals than images. In \cite{marziliano2006blsignals}, the authors propose a joint model of high frequency noise and bandlimited signal to perform exact reconstruction based on the theory of finite rate of innovations, with application to compression of ECG.
    The separation of bioacoustic signals through their spectrogram is studied in \cite{kreme2024separation}, relying on a composite model of non-stationary components.

    A classical strategy for solving ill-posed linear inverse problems is to reconstruct the signal as the solution of an optimization problem. We consider composite optimization problems of the form
    \begin{equation} \label{eq:compositepb}
        (\hfb,\hfh) \in \underset{(s_1,s_2)}{\arg\min} \ \mathcal{D}(\bm{y}, \bm{\Phi}_1(s_1) + \bm{\Phi}_2(s_2)) + \mathcal{R}_1 (s_1) + \mathcal{R}_2(s_2),
    \end{equation}
    in which the data-fidelity $\mathcal{D}(\bm{y}, \bm{\Phi}_1(s_1) + \bm{\Phi}_2(s_2))$ constrains the solution to correspond to the observations and the regularizations $\mathcal{R}_1 (s_1)$ and $\mathcal{R}_2(s_2)$ promote specific behaviors according to some prior knowledge (see Section~\ref{sec:regularizationblabla}). We specifically consider a sparse-plus-smooth model where $\mathcal{R}_1$ is a sparsity-promoting regularization while $\mathcal{R}_2$ favors smooth components.

    We make the assumption that the noise corrupting the data is Gaussian, which motivates the choice of a quadratic data fidelity of the form
    \begin{equation*}
        \mathcal{D} (\bm{y}, \cdot\,) =  \frac{1}{2} \| \bm{y} - \, \cdot \ \|_2^2.
    \end{equation*}
    This choice is essential to our analysis, enabling us to later identify a quadratic optimization subproblem. 
    
    \subsection{From sparse-versus-smooth to sparse-plus-smooth regularization}
    \label{sec:regularizationblabla}

    Consider for a moment the classical single-component version of the optimization problem~\eqref{eq:compositepb}, i.e.,
    \begin{equation} \label{eq:singletpb}
        \underset{s}{\arg\min} \ \mathcal{D}(\bm{y}, \bm{\Phi}(s)) + \mathcal{R} (s)
    \end{equation}
    which solves a non-composite inverse problem. In this work, we distinguish two types of regularization for ill-posed inverse problems.

    Smooth regularization corresponds to the case where the penalty is a quadratic norm over some Hilbert space $\mathcal{H}$, of the typical form 
    \begin{equation*}
        \mathcal{R} (s) = \lambda \| s \|_{\mathcal{H}}^2.
    \end{equation*}
    This type of regularization has been widely used over more than 70 years, with countless applications. The name stems from the smoothing effect it has on the solution, with pioneering works presented in \cite{Kimeldorf1970}
    on the connections between quadratic regularization and spline smoothing. 
    One of the earliest applied versions of quadratic regularization is Tikhonov regularization~\cite{tikhonov1963solution}, also termed ridge-regression in statistics~\cite{Hoerl1962ridge}.  
    Later on, representer theorems have been introduced to precisely specify the form of the solutions, drawing the connection between Hilbertian penalization and the RKHS theory~\cite{Scholkopf2001generalized}. Interestingly, the quadratically-penalized problems are fairly simple to solve as they can be recast into a finite-dimensional formulation, even when the signal $s$ to recover is continuously defined, as for instance with $\mathcal{H} = \ldx$ the space of square integrable functions over a domain $\cX$. 

    The second approach considered in this work is that of sparsity-promoting regularizations
    \begin{equation*}
        \mathcal{R}(s) = \lambda \mathcal{R}^\mathrm{S}(s).
    \end{equation*}
    In this setting, we consider a theoretical model in which the signal $s$ belongs to a Banach search space $\cB$, and the regularization $\cRS: \cB \to \R_+$ is a convex function. This framework is more general than the Hilbertian case discussed previously and encompasses most of the commonly-used sparsity-promoting penalties.
    Doing so, we fit into the formalism presented in \cite{boyer2019representer} and therefore we rely on its abstract representer theorem, which links the choice of the convex penalty $\cRS$ to the structure of the solutions, thereby providing a theoretical explanation for the sparsity-inducing behavior of certain regularizations.
    A central class of examples within this framework is given by Banach norms. When $\mathcal{R}^\mathrm{S} (s) = \| s \|_{\mathcal{B}}$, the solution set of the single-component problem \eqref{eq:singletpb} admits a precise characterization: solutions can be expressed as finite combinations of extreme points of the unit ball of $\cB$ \cite{Unser2020,bredies2020sparsity}. In finite-dimensional Euclidean spaces, the prototypical example is the $\ell_1$-norm, whose extreme points are the canonical basis vectors \cite{tibshirani1996regression,chen2001atomic}. More generally, sparse dictionary reconstruction can be achieved using atomic norms \cite{chandrasekaran2012convex}.
    The Banach-norm perspective further allows one to address continuous-domain sparse reconstruction, for instance performing spikes recovery using the total-variation norm on measures \cite{laville2021sparse} or with the design of specific norms for more complex functional settings \cite{laville2023a,ambrosio2023,decastro2024}.
    However, many classical sparse recovery examples also fit into the $\cRS$-framework without being Banach norms. Positivity constraints, for instance, are known to promote sparse solutions in certain inverse problems \cite{Donoho_Tanner_2005,Slawski_Hein_2011}. Another important model arises from precomposing a sparsity-promoting norm with a differential operator, leading to regularizers that are only seminorms. This includes analysis-prior models \cite{elad2006analysis} and spline reconstruction problems \cite{unser2017splines}.
    Generally, solving sparsity-promoting penalized problems is more challenging than quadratic regularization in Hilbert spaces: there may exist infinitely many solutions in the space $\cB$, and the lack of an inner product structure precludes the use of standard gradient-based methods.

    In this work, we consider the combination of the two above-mentioned penalties, which we refer to as \emph{sparse-plus-smooth} regularization. We focus on problem~\eqref{eq:compositepb} with
        \begin{equation*}
        \mathcal{R}_1 (s_1) = \lambda_1 \mathcal{R}^\mathrm{S}(s_1) \qquad\text{and}\qquad \mathcal{R}_2 (s_2) = \lambda_2 \| s_2 \|_{\mathcal{H}}^2
    \end{equation*}
    for $\cRS$ a general convex sparsity-promoting function $\cB \to \R_+$ and $\lambda_1, \lambda_2 > 0$. This problem has already been studied in \cite{debarre2021continuous} using a Banach norm for $\cRS$, covering continuous-domain reconstruction and regularization operators. The authors demonstrated that the solutions of the composite optimization problem indeed behave as expected: the sparse component $s_1$ admits a sparse structure as if it solved non-composite problem, and similarly the component $s_2$ behaves as the solution of a smooth problem. This result had been exemplified earlier in various finite-dimensional application cases, using sparse-plus-smooth optimization problems with discrete vector components, see \cite{mol2004,gholami2013balanced,debarnot2020learning}. In a discrete setting, it is possible to decouple the resolution of the composite optimization problem: the sparse component can be directly identified with a single-component optimization problem and the smooth component is unique and deduced afterward. The requirements for such a decoupling to take place are presented in \cite{jarret2024decoupled}, considering the general case in which operators are involved within the regularization terms.


    \subsection{Contributions and outline}

    Our main contribution is to identify a new way of conceiving the optimization problem \eqref{eq:compositepb} by decoupling its resolution, extending the result of \cite{jarret2024decoupled} to general Banach and Hilbert spaces. More precisely, we show that the composite optimization problem~\eqref{eq:compositepb} can be reduced to two subproblems. The first one is brought back to the Banach space and directly identifies the sparse component, while the second one is a quadratic optimization problem that admits an explicit solution for the smooth component. We obtain a more precise composite representer theorem than existing results and identify in particular the smooth component in closed form as a function of the vector of measurements of the sparse component.


    We illustrate the decoupling mechanism on a composite deconvolution problem inspired by microscopy imaging.
    We demonstrate in simulations the relevance of using a composite model over a single-component approach for problems involving a smooth background. We also show a significant temporal gain over algorithms exploiting the decoupling of the initial problem compared to the standard approach based on the representer theorem of \cite{debarre2021continuous}. 




\section{Representer Theorems for Single-Component Problems}

    As they will be useful later on, we recall here the classical representer theorems that characterize the solutions of single-component optimization problems such as \eqref{eq:singletpb}, along with the relevant topological structures on the involved search spaces. In particular, smooth regularization can be formulated over regular Hilbert spaces while
    sparse regularization requires to introduce a weak-* topology on Banach spaces (using a dual norm).

    
    \subsection{Smooth regularization of inverse problems}

    {We assume the search space $\mathcal{H}$ to be a Hilbert space with Hilbertian norm $s \mapsto \| s \|_\mathcal{H} := \sqrt{\langle s, s \rangle_{\mathcal{H}}}$. By the Riesz representer theorem, any continuous linear functional $\nu: \cH \to \R$ can be represented as the inner product with an element $\phi_\nu \in \cH$ such that $\nu(f) = \langle \phi_\nu, f \rangle_\cH$ for any $f \in \cH$. The topological dual $\cH'$ is then identified as the search space $\cH$ itself.}

    {
    For $\phihsingle= (\phihh_1 ,\ldots , \phihh_L) \in \mathcal{H}^L$, we denote
    \begin{equation}
        \phihsingle(s) := \left( \langle \phihh_1, s \rangle_{\mathcal{H}} , \ldots , \langle \phihh_L, s \rangle_{\mathcal{H}} \right). 
        \label{eq:phi-hilbert}
    \end{equation}
    The adjoint of $\phihsingle$ is the only operator $\phihsingle^* : \R^L \rightarrow \mathcal{H}$ which satisfies the equality 
    \begin{equation*}
        \langle \phih(s), \bm{y}\rangle_{\R^L} = \langle s, \phih^*(\bm{y})\rangle_\cH,
    \end{equation*}
    for any $s \in \cH$ and $\bm{y}\in\R^L$.
    Its expression is given by $\phihsingle^* (\bm{y}) = \sum_{1 \leq \ell \leq L} y_\ell \phihh_\ell$.
    The entries of the Gram matrix $\mathrm{G} := \phihsingle \phihsingle^* \in \R^{L\times L}$ are given by $G[k,\ell] = \langle \phihh_k , \phihh_\ell \rangle_\cH$.}

    
   Proposition \ref{prop:hilbertRT} hereafter is a classical result for quadratic optimization over Hilbert spaces. Significantly more general formulations can be found in~\cite{Scholkopf2001generalized} but we restrict here to the case which is relevant for our purpose. We follow the exposition of~\cite[Theorem~7]{caponera2021nonparametric}, which is based on \cite[Section~3.2]{Unser2020}. 

    \begin{proposition}[Representer theorem on Hilbert spaces]
    \label{prop:hilbertRT}
    Let $\bm{y}\in \R^L$, $\phihsingle= (\phihh_1 ,\ldots , \phihh_L) \in \mathcal{H}^L$, and $\lambda > 0$. If the measurement functionals $\phihh_1, \dots, \phihh_L$ are linearly independent, the optimization problem 
    \begin{equation} \label{eq:hibertpb}
        \inf_{s \in \mathcal{H}} \frac{1}{2} \| \bm{y} - \phihsingle(s) \|_2^2 + \frac{\lambda}{2} \| s \|_{\mathcal{H}}^2
    \end{equation}
    admits a unique solution $\hfh \in \mathcal{H}$ which is given by
    \begin{equation} \label{eq:widehathilb}
        \hfh  = \phihsingle^* ( \phihsingle \phihsingle^* + \lambda \mathbf{I}_L)^{-1} \bm{y}.
    \end{equation}
    \end{proposition}
    
    \noindent In other terms, Proposition~\ref{prop:hilbertRT} states that the solution set of ~\eqref{eq:hibertpb} is 
    \begin{equation}
        \mathcal{U} (\lambda) := \{ \phihsingle^* ( \phihsingle \phihsingle^* + \lambda \mathbf{I}_L)^{-1} \bm{y}\}.
    \end{equation}
    The existence and uniqueness of the solution of \eqref{eq:hibertpb} follows from the the fact that the functional of \eqref{eq:hibertpb} is continuous and strongly convex due to the quadratic regularization.
    The vector $ \hfh $ is in the span of the functionals $\phihh_\ell$ since \eqref{eq:widehathilb} can be reinterpreted as 
    \begin{equation*}
    \hfh  = \sum_{1 \leq \ell \leq L} \alpha_\ell \phihh_\ell \quad \text{where} \quad \bm{\alpha} = ( \phihsingle \phihsingle^* + \lambda \mathbf{I}_L)^{-1} \bm{y} \in \R^L. 
    \end{equation*}


    \subsection{Sparse regularization of inverse problems}

    The solution resulting from the quadratic Hilbertian regularization is generally not sparse (it can be interpreted as an orthogonal projection over a Hilbert norm ball, see \cite{engl1996regularization}).
    A general framework for sparse reconstruction is classically formalized using a non-reflexive Banach space $\mathcal{B}$ as the search space, for which we assume the existence of a predual space \cite{unser2021unifying}.
    
    We fix two Banach spaces $( \mathcal{A}, \| \cdot \|_{\mathcal{A}})$ and $( \mathcal{B}, \| \cdot \|_{\mathcal{B}})$ such that $\mathcal{B} = \mathcal{A}'$ is the topological dual of $\mathcal{A}$ and $\| \cdot \|_{\mathcal{B}}$ is the dual norm 
    \begin{equation*}
        \| s \|_{\mathcal{B}} = \sup_{\|u \|_{\mathcal{A}}=1} \langle u, s \rangle_{\mathcal{A}\times\mathcal{B}}, 
    \end{equation*}
    {where the duality product is denoted as $\langle u, s \rangle_{\mathcal{A}\times\mathcal{B}} = s(u)$.}
    This assumption implies that $\mathcal{B}$ can be endowed with the weak*-topology inherited from $\mathcal{A}$. We say that the sequence $(s_n)_{n\geq 1}$ of elements $s_n \in \mathcal{B}$ converges to $s \in \mathcal{B}$ for the weak*-topology if 
    \begin{equation*}
        \langle u , s_n \rangle_{\mathcal{A}\times\mathcal{B}} \underset{n\rightarrow \infty}{\longrightarrow} \langle u, s \rangle_{\mathcal{A}\times\mathcal{B}} 
    \end{equation*}
    for any $u \in \mathcal{A}$.
    In other terms, $\mathcal{B}$ admits a predual\footnotemark.
    \footnotetext{The Hilbert scenario presented in the previous section corresponds to the case $\mathcal{A} = \mathcal{B} = \mathcal{H}$.}

    A typical example is the space of Radon measures $(\mathcal{B}, \|\cdot \|_{\mathcal{B}}) = (\mx, \|\cdot \|_{\mathcal{M}})$ over a continuous domain $\cX$. Its predual is the Banach space of continuous vanishing functions for the supremum norm $(\mathcal{A}, \|\cdot \|_{\mathcal{A}}) =(\czx, \| \cdot \|_{\infty})$. The Radon measures remarkably contain Dirac impulses, which make them well-suited for continuous-domain sparse recovery. For any $d\in\mathbb{N}^*$, the case $\cX=\mathbb{R}^d$ is for instance presented in~\cite{unser2017splines,Unser2020} and the periodic case $\cX=\mathbb{T}^d$ is covered by \cite{fageot2020tv}.

    Without the Hilbert space structure, the measurement operator $\phibsingle= (\phibb_1, \ldots, \phibb_L)$ is made of measurement functionals $\phibb_\ell \in \mathcal{A}$ from the predual space. Then, $\phibsingle$ specifies a linear and weak*-continuous mapping $\phibsingle: \mathcal{A}' =\mathcal{B} \rightarrow \R^L$ via 
    \begin{equation}
        \label{eq:phi-banach}
        \phibsingle (s) := (\langle \phibb_\ell, s \rangle_{\mathcal{A}\times\mathcal{B}})_{1\leq \ell \leq L} \in \mathbb{R}^L.
    \end{equation}
    This definition ensures that $\phibb$ is a bounded operator $\left(\cB, \norm{\cdot}_\cB \right) \to \left(\R^L, \norm{\cdot}_2\right)$, as for any $s\in\cB$ the duality inequality holds: $\lvert \langle \phibb_\ell, s \rangle_{\mathcal{A}\times\mathcal{B}} \rvert \leq \norm{\phibb_\ell}_\cA \norm{s}_\cB$.

    Within this framework, the solutions of the single-component optimization problem~\eqref{eq:singletpb} with $s\in\cB$  are significantly influenced by the choice of the sparsity-promoting penalty $\mathcal{R}(s) = \lambda \cRS(s)$. Remarkably, we have access to abstract representer theorems to characterize the solution set. We follow the presentation of \cite{boyer2019representer} but similar results were proposed at the same time with somewhat different proofs in \cite{bredies2020sparsity}. The abstract theorem presented hereafter generalizes Proposition~\ref{prop:hilbertRT} on Hilbert spaces.
    
    \begin{proposition}[Abstract representer theorem for penalized problems (Corollary 2 in \cite{boyer2019representer})]
    \label{prop:abstract-rt}
    Let $\cRS: \cB \to \R_+$ be a lower semi-continuous convex function, $\bm{y}\in \R^L$, $\phibsingle= (\phibb_1 ,\ldots , \phibb_L) \in \mathcal{A}^L$, and $\lambda > 0$. Moreover, let us assume that the optimization problem
    \begin{equation} \label{eq:sparse-rt}
        \inf_{s \in \mathcal{B}}\ \frac{1}{2} \| \bm{y} - \phibsingle(s) \|_2^2 + \lambda \cRS(s)
    \end{equation}
    has a nonempty solution set $\mathcal{V}(\lambda)$.
    \begin{itemize}
        \item First, any solution $\widehat{s} \in \mathcal{V}(\lambda)$ shares the same measurement vector, i.e., $\phib \left(\mathcal{V}(\lambda)\right)$ is a singleton.
        \item Second, any solution $\widehat{s} \in \mathcal{V}(\lambda)$ leads to the same value of the penalty $R^\star := \cRS(\widehat{s})$.
        \item Third, the extreme points of $\mathcal{V}(\lambda)$ can be expressed as convex combinations of at most $L$ extreme points of the level set $C^\star:= \left\{s \in \cB,\ \cRS(s) \leq R^\star\right\}$ (or $L+1$ in the specific case $R^\star = \inf_\cB \cRS$). 
    \end{itemize}
    \end{proposition}

    \begin{proof}
        The first statement repeats Remark~6 in \cite{boyer2019representer}. The proof relies on the strict convexity of the data-fidelity term, as presented in Lemma~1 from \cite{tibshirani2013lasso}. The second statement is the consequence of the first one: two different solutions $\widehat{r}, \widehat{s} \in \mathcal{V}(\lambda)$ have the same fitted value $\phib(\widehat{r}) = \phib(\widehat{s})$, hence they also produce the same value of the quadratic data-fidelity term. As they are both solutions to the problem, they reach the same global minimum and so they necessarily lead to the same value of the penalty term $\cRS(\widehat{r}) = \cRS(\widehat{s})$, given that $\lambda > 0$. The third statement is the direct application of Corollary~2 from \cite{boyer2019representer}, using the lower semi-continuity of $\cRS$ to ensure that $C^\star$ is linearly closed (Remark~5 in the article). 
    \end{proof}

    Proposition~\ref{prop:abstract-rt} sheds light on the mechanics behind sparsity-promoting penalties: sparse reconstruction is naturally obtained using a function $\cR^{\mathrm{S}}$ whose level sets have sparse extreme points, and the reconstruction atoms can be chosen to be the extreme points themselves.

    More precise representer theorems can be deduced from Proposition~\ref{prop:abstract-rt}, for instance considering (generalized) total-variation norms for measures and splines reconstructions\cite{Fisher1975,Unser2016representer,unser2017splines,bredies2020sparsity,debarre2022sparsest}. As an example, we recall the specific case of a Banach-norm regularization $\cRS = \norm{s}_\cB$ in Proposition~\ref{prop:banachRT}, which combines existence and topological results from~\cite[Proposition 8]{gupta2018continuous} and the extreme point characterization of \cite[Theorem 3.1]{boyer2019representer}.

    \begin{proposition}[Representer theorem on Banach spaces]
    \label{prop:banachRT}
    Let $\bm{y}\in \R^L$, $\phibsingle= (\phibb_1 ,\ldots , \phibb_L) \in \mathcal{A}^L$, and $\lambda > 0$. The optimization problem 
    \begin{equation} \label{eq:banachpb}
        \inf_{s \in \mathcal{B}} \frac{1}{2} \| \bm{y} - \phibsingle(s) \|_2^2 + \lambda \| s \|_{\mathcal{B}}
    \end{equation}
    admits at least a solution and its solution set $\mathcal{V} (\lambda)$ is weak*-compact, convex, and is the closed convex hull of its extreme points.

    Any extreme point $\hfb$ of $\mathcal{V} (\lambda)$ is such that 
    \begin{equation} \label{eq:widehatbanachextreme}
        \hfb  = \sum_{1 \leq k \leq K} \alpha_k e_k
    \end{equation}
    where $e_k$ are distinct extreme points of the unit ball of $\|\cdot \|_{\mathcal{B}}$,  $\alpha_k \in \mathbb{R}$, and $0 \leq K \leq L$.
    \end{proposition}

\section{Sparse-plus-Smooth Composite Representer Theorem}

    Building on the analysis of the single-component problems, we now turn to the theoretical study of the composite case.

    \subsection{Main theorem}

    Let us recall the composite problem \eqref{eq:compositepb} with a quadratic data-fidelity term. For $\lambda_1, \lambda_2 > 0$, the objective functional is given by
    \begin{equation}
        \label{eq:composite-optim}
        \mathcal{J}(\fb, \fh) := \frac{1}{2} \| \bm{y} - (\phib (\fb) + \phih (\fh)) \|_2^2  + \lb \cRS(\fb) + \frac{\lh}{2} \| \fh \|_{\mathcal{H}}^2,
    \end{equation}
    where $\phib \in \mathcal{A}^L$ and $\phih\in\cH^L$ respectively sample the sparse and the smooth components. 
    The set of pairs of minimizers is defined as
    \begin{equation} \label{eq:argminisback}
        \mathcal{W} (\lb,\lh) := \underset{(\fb,\fh) \in \mathcal{B} \times \mathcal{H}}{\arg\min} \mathcal{J}(\fb, \fh).
    \end{equation}
    We need to define the matrix
    \begin{equation}
        \label{eq:def-Mphi}
        \mathbf{M}_{\lh} := \frac{1}{\lh} \left(\phih\phih^* + \lh \mathbf{I}_L\right) = \frac{1}{\lh} \left( \langle \phihh_k , \phihh_\ell \rangle_{\mathcal{H}} + \lh \delta[k - \ell] \right)_{1\leq k, \ell \leq L}.
    \end{equation}
    This matrix is positive definite. It is therefore invertible and admits a square-root matrix.

    Our main result, Theorem~\ref{theo:main} hereafter, reduces the analysis of a composite optimization problem over $\mathcal{B}\times\mathcal{H}$ to a problem over $\mathcal{B}$, hence decoupling the contributions of the two components. The proof is given in Appendix~\ref{app:prooftheo1}.

    \begin{theorem} \label{theo:main}
        Let $\cRS: \cB \to \R_+$ be a convex function,  $\bm{y}\in \R^L$, $\lambda_1, \lh > 0$, $\phib \in \mathcal{A}^L$ and $\phih \in \mathcal{H}^L$.
        The solution set $\mathcal{W} (\lb,\lh)$ can be written as
        \begin{equation}
            \label{eq:decoupling}
            \mathcal{W} (\lb,\lh) = \mathcal{V}(\mathbf{M}_{\lh},\lb) \times \{\hfh\}
        \end{equation}
        with
        \begin{align}
            \mathcal{V}(\mathbf{M}_{\lh},\lb) &:= \underset{ \fb \in \mathcal{B}}{\arg\min} \quad \|\mathbf{M}_{\lh}^{-\frac{1}{2}} ( \bm{y} - \phib(\fb)  ) \|_2^2  + \lb \cRS(\fb), \label{eq:banachpart} \\
            \hfh &:=  \frac{1}{\lh}\phih^* \mathbf{M}_{\lh}^{-1} \left(\bm{y} - \bm{w}\right), \label{eq:hilbertpart}
        \end{align}
        where the vector $\bm{w} := \phib ( \hfb)$ is unique and independent of the solution $\hfb \in \mathcal{V}(\mathbf{M}_{\lh},\lb)$. 
    \end{theorem}

    Theorem~\ref{theo:main} reveals the interests of a composite framework using sparsity-promoting and quadratic Hilbertian regularizations. The sparse component is made of extreme points of the level sets of $\cRS$ while the smooth component lives in the finite-dimensional space of the measurement functionals $(\phihh_1, \ldots, \phihh_L)$. It implies in particular that the general form of an extreme point solution is 
    \begin{equation*}
        \left(\hfb, \hfh \right) = \left(\sum_{k} \alpha_k e_k, \sum_{1 \leq \ell \leq L} \beta_\ell \phihh_{\ell}\right),
    \end{equation*}
    where the $e_k$ are distinct extreme points of the the level sets of $\cRS$. 


    Another consequence of Theorem~\ref{theo:main} is that known results over $\mathcal{B}$ are directly transferred into $\mathcal{B}\times \mathcal{H}$, as the decoupling holds for any convex penalty $\cRS$. The properties of the sparse component directly stem from the shape of the solutions in $\mathcal{V}(\mathbf{M}_{\lh}, \lb)$. We illustrate this principle with a Banach-norm regularization in the following corollary, stating a representer theorem which characterizes the extreme point solutions.

    \begin{corollary}[Decoupling representer theorem]
    \label{cor:main}
    Let us consider the problem~\eqref{eq:argminisback} with $\cRS(\cdot) = \norm{\cdot}_\cB $ over $\cB$. Invoking Proposition~\ref{prop:banachRT} to $\mathcal{V}(\mathbf{M}_{\lh},\lb)$ in \eqref{eq:banachpart}, the solution set $\mathcal{W} (\lb,\lh)$ is non-empty, convex, and weak*-compact in $\mathcal{B}\times \mathcal{H}$. Its extreme points are of the form 
    \begin{equation} \label{eq:extreme}
        (\hfb , \hfh ) = \left( \sum_{1 \leq k \leq K} \alpha_k e_k , \frac{1}{\lh}\phih^* \mathbf{M}_{\lh}^{-1} \left(\bm{y} - \bm{w}\right)\right)
    \end{equation}
    where $\alpha_k \neq 0$, $e_k$ are distinct extreme points of the unit ball $\{ \fb \in \mathcal{B}, \ \| \fb \|_{\mathcal{B}} \leq 1\}$, $0 \leq K \leq L$, and $\hfh$ is given by \eqref{eq:hilbertpart}. 
    \end{corollary}

    Equation~\eqref{eq:extreme} was already a consequence of~\cite[Theorem 2]{unser2022convex}, but Theorem~\ref{theo:main} specifies the exact form of the smooth component $\hfb$. Incidentally, spline-based reconstruction could be readily considered for the sparse component using $\cRS(\cdot) = \|\mathrm{L}\cdot\|_\mathcal{B}$, with $\mathrm{L}$ a pseudo-differential operator \cite{unser2017splines}.

    \subsection{Maximum regularization parameter for TV-norm}

    Setting the regularization parameters in optimization problems is a notoriously sensitive task. 
    If we consider a single-component problem with a Hilbert penalty, there exists an optimal value of the regularization parameter when the source signal follows a random Gaussian model \cite{Badou}. Beyond this ideal model, finding a relevant value for this parameter is still an open question and many strategies have been proposed in the literature \cite{hansen2000,park2008parameter}. The case of single-component Banach problems is not simpler \cite{deladalle2014sugar,chirinos2024parameter}. 
    However, when an $\ell_1$-norm or a total-variation norm is considered, the relevant values of the regularization parameter $\lambda > 0$ are confined to an explicit interval. Indeed, there exists a problem-dependent maximum value $\lambda_{\mathrm{max}} > 0$ above which the solution of the optimization is unique and reduced to the null signal. This is a well-known result for LASSO-type problems, both in discrete \cite{tibshirani2013lasso}, \cite[Proposition~II.1]{koulouri2021} and continuous settings \cite[Proposition~10]{debarre2022sparsest}, \cite[Proposition~4.3]{jarret_grid_2025}.

    We first extend this result to the generic case of single-components problems penalized with a Banach norm in the following theorem, whose proof is deferred to Appendix~\ref{app:proof:2}.
    \begin{theorem}[Regularization bound on Banach-penalized problems]
    \label{prop:single-lmax}
    Let $\bm{y}\in \R^L$ and $\phib \in \cA^L$.
    For $\lambda >0$, we consider the solution set of the single-component optimization problem
    \begin{equation}
        \label{eq:single-banach}
        \mathcal{V}(\lambda) := \underset{s \in \cB}{\arg \min} \frac{1}{2} \norm{\bm{y} - \phib(s)}_2^2 + \lambda \norm{s}_{\cB}.
    \end{equation}
    Let us define
    \begin{equation}
        \lambda_{\mathrm{max}} = \norm{\phib^{*}\bm{y}}_{\mathcal{A}}.
    \end{equation}
    The two statements hold:
    \begin{itemize}
        \item[1.] \label{item:lm:1} For any $\lambda \geq \lm$, the solution to \eqref{eq:single-banach} is unique and the solution set is reduced to $\mathcal{V}(\lambda) = \left\{0\right\} \subset \cB$.
        \item[2.] \label{item:lm:2} For any $\lambda < \lm$, the null element $0$ is not solution, i.e., $0 \notin \mathcal{V}(\lambda)$.
    \end{itemize}
    \end{theorem}

    This result on the maximum value of the regularization parameter can be transferred to composite problems using Theorem~\ref{theo:main}. It reveals an explicit dependence between the two regularization parameters of composite problems, as illustrated in Proposition~\ref{prop:lmax}.

    \begin{proposition}[Maximum value of $\lambda_1$]
    \label{prop:lmax}
    Consider the composite optimization problem \eqref{eq:argminisback} with $\cRS(\cdot) = \lVert \cdot \rVert_\mathcal{B}$. For a fixed $\lh >0$, we define
    \begin{equation}
        \label{eq:l1max}
        \lambda_{1, \mathrm{max}} := \lVert \phib^* \mathbf{M}_{\lambda_2}^{-1} \bm{y} \rVert_\cA = \lambda_2 \lVert \phib^* \left( \phih \phih^* + \lambda_2 \mathbf{I}_L \right)^{-1} \bm{y} \rVert_\cA.
    \end{equation}
    For any $\lambda_1 \geq \lambda_{1, \mathrm{max}}$, the solution set for the Banach component is reduced to the singleton zero 
    $$\mathcal{V}(\mathbf{M}_{\lh}, \lb) = \left\{ 0 \right\}$$
    and Problem \eqref{eq:argminisback} is equivalent to a single-component Hilbert problem.
    \end{proposition}

    \begin{proof}
        The sparse component subproblem~\eqref{eq:banachpart} is a Banach-norm penalized single-component problem as in~\eqref{eq:single-banach} with operator $\mathbf{M}_{\lh}^{-\frac{1}{2}} \phib$ and measurement vector $\mathbf{M}_{\lh}^{-\frac{1}{2}}\bm{y}$, hence the maximum value of the optimization parameter writes as
        \begin{equation*}
            \lambda_{1, \mathrm{max}} = \lVert (\mathbf{M}_{\lh}^{-\frac{1}{2}} \phib)^*  (\mathbf{M}_{\lh}^{-\frac{1}{2}} \bm{y}) \rVert_\cA = \lVert \phib^* (\mathbf{M}_{\lh}^{-\frac{1}{2}})^* \mathbf{M}_{\lh}^{-\frac{1}{2}} \bm{y} \rVert_\cA = \lVert \phib^* \mathbf{M}_{\lambda_2}^{-1} \bm{y} \rVert_\cA,
        \end{equation*}
        using the symmetry of the matrix $\mathbf{M}_{\lambda_2}$.
    \end{proof}

    Equation \eqref{eq:l1max} demonstrates a natural dependence of $\lb$ onto $\lh$. Indeed, the range of relevant values for $\lb$ is $[0, \lambda_{1, \mathrm{max}}]$ which depends on $\lh$. Hence, by first choosing $\lh$ then setting $\lb = \alpha \lambda_{1, \mathrm{max}}$ for $0 < \alpha < 1$ we ensure that the value of the regularization parameters is consistent with the problem. This approach was already proposed and discussed in our previous work on composite problems \cite{jarret2024decoupled}. We illustrate this dependency and how our scaling rule allows to decouple the choice of parameters with a cross table of reconstructions in Appendix~\ref{app:reconstructions}.


    Additionally, equation \eqref{eq:l1max} provides information on the asymptotic behavior of $\lbm$. When $\lambda_2$ is small, $\lh \mathbf{I}_L$ is negligible compared to $\phih \phih^*$ and so $\lbm$ is proportional to $\lh$. However, for large values of $\lh$, the Gram matrix $\phih \phih^*$ is itself negligible and $\lbm$ tends to $\lVert \phib^* \bm{y}\rVert_\cA$. This latter result is consistent with the maximum value of the regularization parameter for the single-component Banach-penalized problem~\eqref{eq:single-banach}.


\section{Application to Super-Resolved Deconvolution}
\label{sec:application}

    Deconvolution is a classical imaging inverse problem, appearing in various practical applications such as microscopy, photography or astronomy. It consists in recovering a high resolution image from blurred observations, possibly corrupted by various artifacts (measurement noise, external illumination, etc).
    The presence of background information in the signal to recover, even with a small intensity, significantly deteriorates the quality of traditional sparsity-based single-component reconstruction methods.

    In this section, we illustrate our representer theorem for a deconvolution optimization problem. We first emphasize the decoupling in the specific case of \emph{measures-plus-Hilbert} signal before performing actual reconstructions using a grid-based approximation of continuous-domain sparse signals.


    \subsection{Definition of the composite model}
    \label{sec:comp:model}
        

        We consider a scenario inspired from microscopy imaging, in which the measurements of the sample may be corrupted by the presence of out-of-focus elements. Such a situation can be accurately modeled using a composite model and in particular a smooth background component.
        We define a simplified deconvolution inverse problem, in which the measurements $\bm{y}\in\R^L$ are expressed as
        \begin{equation*}
            \bm{y} \approx S \left\{g * (s_1^\dagger, s_2^\dagger)\right\}.
        \end{equation*}
        The ground truth signal $s^\dagger := (s_1^\dagger, s_2^\dagger)$ is composed of a sparse foreground component $s_1^\dagger \in \mathcal{B}$ and a smooth background component $s_2^\dagger \in \mathcal{H}$. The action of the microscope on the scene is modeled as the convolution with the \emph{point-spread function} $g$ (PSF), before sampling on the grid of observation pixels with the operator $S$.

        Remarkably, the signals to recover can be chosen as continuously-defined functions, a natural way to perform super-resolution recovery. A motivation for such a continuous model comes for instance from Single-Molecule Localization Microscopy (SMLM) \cite{sage2019smlm}, a microscope modality specifically known for achieving reconstruction below the diffraction limit \cite{huang2017super,denoyelle2019sliding,laville2021sparse}.


        \subsubsection{Function spaces}
            For $d \in \mathbb{N}^*$, let $\cX = \mathbb{R}^d$ be the domain of the signal to recover $s^\dagger$. 

            We consider the classical model of Radon measures $\mathcal{B} = \mx$, in which sparse signals can be represented with Dirac impulses. The {total-variation norm} $\norm{\cdot}_\mathcal{M}$ on $\mx$ is known to promote sparse solutions when used as a regularizer for optimization problems \cite{unser2017splines}. For $\czx$ the space of continuous vanishing functions on $\cX$, which is the Banach predual of $\mx$, we recall the dual definition of the total-variation norm
            \begin{equation*}
                \forall m \in \mx, \quad \lVert m \rVert_\mathcal{M} = \underset{\substack{\varphi \in \czx \\ \| \varphi \|_\infty = 1}}{\sup} \int_\cX \varphi(x) \mathrm{d}m(x).
            \end{equation*}
            

            For the smooth background component $s_2$, we consider the Hilbert space $\cH = \cHk$ the RKHS\footnote{\emph{Reproducing Kernel Hilbert Space}.} induced by the Gaussian kernel $k: \cX \times \cX \to \R$, also called \emph{Gaussian radial basis function}, defined as
            \begin{equation*}
                \forall x, y \in \cX, \qquad k(x, y) = g_0(\norm{y - x}) := \frac{1}{(\sqrt{2\pi\sigma_0^2})^d} \exp{ \left({-\frac{\norm{y - x}^2}{2 \sigma_0^2}}\right)},
            \end{equation*}
            for $\sigma_0 > 0$. The kernel $k$ defines an inner product $\langle \cdot, \cdot \rangle_{\cHk}$ and the associated Hilbert norm $\norm{\cdot}_{\cHk}$ on $\cHk$. The formal construction of the RKHS is technical and we refer for instance to \cite[Chapter~10]{wendland2004}. Notably, no direct expression for the inner product  $\langle f, g \rangle_{\cHk}$ is available in the general case of arbitrary $f, g \in \cHk$. This Gaussian RKHS $\cHk$ is a subspace of the square integrable functions $\ldx$, containing the functions with spectral decay dominated by the decay of $k$ (see, e.g., \cite{minh2010} for a more detailed characterization of $\cHk$). We later rely on this property to adjust the intensity of the promoted smoothness and tune the reconstruction with $\sigma_0$.

            Additionally, we assume that the ground truth signals $s_1^\dagger$ and $ s_2^\dagger$ have a compact support, included in the field of view of the microscope. Without loss of generality, we consider this support to be included in the interval $[0, 1]^d$.

        \subsubsection{Measurement operator}        
            In an optical system, each pixel of the sensor collects part of the light coming from the scene and produces one measurement $y_\ell$. This procedure is classically modeled with a convolution between the PSF of the measurement device $g \in \cHk \cap \czx$ and the observed signal $s = (\fb, \fh)$ followed by local integration on a uniform grid \cite{denoyelle2019sliding}. For simplicity, we replace the integration step with a pointwise evaluation of the convolution on the grid.
        
            The sensor is made as a $d$-dimensional uniform grid of $K$ pixels per dimension, leading to $L = K^d$ measurements. The pixel locations are noted as $(x_\ell)_{\ell = 1, \dots, L}$ with $x_\ell \in [0, 1]^d$. For instance, $d=1$ leads to $L=K$ and $x_\ell = {\ell}/({L-1})$ for $\ell \in \{0, \dots, L-1\}$.
        
            It leads to the following expression for the two measurement operators of \eqref{eq:composite-optim}. Consider $\ell \in \{1, \dots, L\}$.
            \begin{itemize}
                \item Definition of $\phib$:
            \begin{equation}
                \label{eq:phib-conv}
                \begin{split}
                \forall \fb \in \mx, \qquad \left[\phib(\fb)\right]_\ell &= \left(g*\fb\right)(x_\ell) \\
                &= \int_\cX g(x_\ell - x) \mathrm{d}\fb(x).
            \end{split}
            \end{equation}
            Keeping the notation of equation \eqref{eq:phi-banach}, this measurement operator corresponds to the measurement functionals
            $$\phibb_\ell = g(x_\ell - \cdot) \in \czx.$$
                \item Definition of $\phih$:
            \begin{equation}
                \label{eq:phih-conv}
                \begin{split}
                \forall \fh \in \cHk, \qquad\left[\phih(\fh)\right]_\ell &= \left(g*\fh\right)(x_\ell) \\
                &= \int_\cX g(x_\ell - x) \fh(x) \mathrm{d}x.
                \end{split}
            \end{equation}
            Expression \eqref{eq:phih-conv}, given in integral form, can also be written using the inner product on $\cHk$ in accordance with \eqref{eq:phi-hilbert}. The identification of the elements $\phihh_\ell$ is deferred to Proposition~\ref{prop:phihh}.
            \end{itemize}

            {
            \begin{proposition}[Lemma~10 in \cite{Berlinet2011reproducing}]
                \label{prop:phihh}
                For any $\fh \in \cHk$, for $\ell \in \left\{ 1, \dots, L \right\}$, the measurement operator can be written as
                \begin{equation}
                    \label{eq:defPhihStar}
                    [\phih(\fh)]_\ell = \langle \phihh_\ell, \fh \rangle_{\cHk}
                \end{equation}
                with $\phihh_\ell$ being the convolution
                \begin{equation*}
                    \phihh_\ell = g_0 * g(x_\ell - \cdot).
                \end{equation*}
            \end{proposition}
            \begin{proof}
                

                By Riesz representer theorem, for every coordinate  $\ell \in \left\{1, \dots, L \right\}$, there exists an element $\phihh_\ell \in \cHk$ such that for any $\fh \in \cHk$ it holds $[\phih(\fh)]_\ell = \langle \phihh_\ell , \fh \rangle_{\cHk}$. In particular, for $t \in \cX$, using $\fh = k(\cdot, t)$ and the reproducing property, we obtain
                \begin{align*}
                    \forall t \in \cX, \qquad \phihh_\ell(t)
                    &= \langle \phihh_\ell, k(\cdot, t) \rangle_{\cHk} = [\phih(k(\cdot, t))]_\ell \\
                    & = \int_\cX g(x_\ell - x) k(x, t) \mathrm{d}x = \int_\cX g(x_\ell - x) g_0(t - x) \mathrm{d}x \\
                    & = (g(x_\ell - \cdot) * g_0)(t).
                \end{align*}   
            \end{proof}
            }
            
            We consider the simplistic model of the PSF being an isotropic Gaussian function of known standard deviation $\sigma > 0$ so that
            \begin{equation}
                \label{eq:psf}
                g(x) = \frac{1}{(\sqrt{2\pi\sigma^2})^d} \exp\left({-\frac{\norm{x}^2}{2 \sigma^2}}\right)
            \end{equation}
            with $\norm{\cdot}$ the Euclidian norm on $\R^d$. It corresponds to a simplified version of the classical astigmatism model as presented in \cite{huang2017super,denoyelle2019sliding}.
            
\begin{figure}
                \centering
                \includegraphics[width=.8\linewidth]{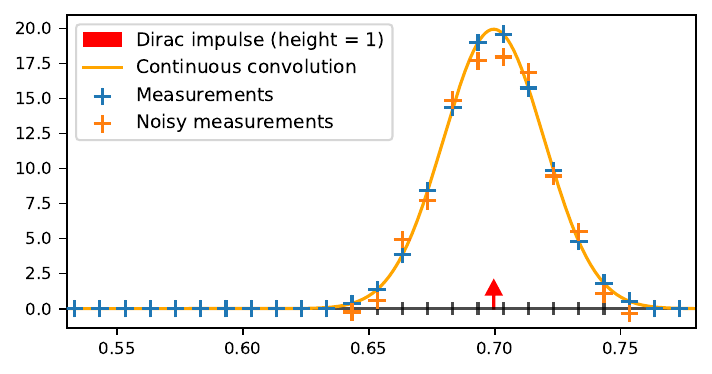}
                \caption{Illustration of the effect of the measurement operator $\phib$ applied to a Dirac signal $\fb^\dagger = \delta_{x_0}$ with $x_0 \approx 0.695$.}
                \label{fig:simple:measop}
            \end{figure}

            Finally, the observations $\bm{y} \in \R^{L}$ are assumed to be corrupted by additive noise as follows
            \begin{equation}
                \label{eq:ip-bis}
                \bm{y} = \phib(s_1^\dagger) + \phih(s_2^\dagger) + \bm{n} \in \R^L
            \end{equation}
            where $\bm{n}$ is a white Gaussian noise of unknown variance $\sigma_{\bm{n}}^2$.
            An illustrative one-dimensional example of simulated measurements is provided in figure~\ref{fig:simple:measop}.



            {
            \begin{remark}[Synthesis operation within $\cHk$]
            In this convolution model, the adjoint operator $\phih^*$ as it appears in \eqref{eq:def-Mphi} and \eqref{eq:hilbertpart} admits the following expression, resulting from the RKHS structure of $\cHk$:
            \begin{equation}
                \label{eq:rkhs-atom}
                \forall \mathbf{h} \in \mathbb{R}^L, \quad \phih^* (\mathbf{h}) = \sum_{\ell} h_\ell\, g_t(x_\ell - \cdot),
            \end{equation}
            with $g_t := g_0 * g$ being another Gaussian kernel of target variance $\sigma_t^2 = \sigma_0^2 + \sigma^2$. In other words, this setup advantageously provides an explicit characterization of the spread of the reconstruction Gaussian kernel, which can be tuned by the practitioner and larger than the one of the measurement kernel.
            \end{remark}
            }

    \subsection{Composite representer theorem}

        As proposed, we address the deconvolution and signal separation task from \eqref{eq:ip-bis} using the two-variables optimization problem \eqref{eq:argminisback}. For a pair of parameters $\lambda_1, \lambda_2 > 0$ we consider

        \begin{equation}
        \label{eq:argmin-sps}
            \underset{(\fb,\fh) \in \mx \times \cHk}{\arg\min} \frac{1}{2} \| \bm{y} - (\phib (\fb) + \phih (\fh)) \|_2^2  + \lb \| \fb \|_{\mathcal{M}} + \frac{\lh}{2} \| \fh \|_{\cHk}^2.
        \end{equation}

        The solutions are indeed decoupled, according to the following result.
        \begin{proposition}
            \label{prop:rt-applied}
            The solution set of problem~\eqref{eq:argmin-sps} is non-empty, convex and weak*-compact.
            Any extreme point solution to the optimization can be expressed as a pair of components $(\hfb, \hfh) \in \mx \times \cHk$ such that
            \begin{align}
                \hfb &= \sum_{k=1}^{K_0} a_k \delta_{z_k}, \label{eq:deconv-banachpart}\\
                \hfh &=  \frac{1}{\lh}\phih^* \mathbf{M}_{\lh}^{-1} \left(\bm{y} - \bm{w}\right),\label{eq:deconv-hilbertpart}
            \end{align}
            in which $(a_k, z_k)_k \in \left(\R \times \cX\right)^{K_0} $ are amplitude-location pairs, $K_0 \leq L$, $\bm{w} = \phib(\hfb)$ is independent on the actual solution $\hfb$, and the matrix $\mathbf{M}_{\lambda_2} \in \R^{L \times L}$ is defined in equation \eqref{eq:def-Mphi}.
        \end{proposition}

        \begin{proof}
            By construction, Corollary~\ref{cor:main} holds for problem \eqref{eq:argmin-sps}. It directly leads to the expression of $\hfh$ in $\eqref{eq:deconv-hilbertpart}$. Regarding the sparse component, we have that
            \begin{equation*}
                \hfb \in \underset{ s_1\in\mx}{\arg\min} \quad  \mathcal{J}_{\mathbf{M}}(s_1)
            \end{equation*}
            with $\mathcal{J}_{\mathbf{M}}(s_1) := \|\mathbf{M}_{\lh}^{-\frac{1}{2}} ( \bm{y} - \phib(s_1)  ) \|_2^2  + \lb \| s_1 \|_{\mathcal{M}}$. From Proposition~\ref{prop:banachRT}, we know that this latter problem admits sparse solutions and that some can be expressed according to \eqref{eq:deconv-banachpart}.
        \end{proof}

    
        
        The smooth background problem \eqref{eq:deconv-hilbertpart} is explicit, so that the only numerical challenge consists in estimating the positions $z_k \in \cX$ and the associated intensities $a_k \in \R$ of the foreground component in \eqref{eq:deconv-banachpart}. The positions are notably more challenging to estimate than the intensities due to their continuously-defined nature and the nonlinear dependence on the measurements. This is done numerically by optimizing the single-component functional $\mathcal{J}_{\mathbf{M}}(\cdot)$ over $\mx$.
        The unique solution background component $\hfh$ is then synthesized after computation of the residuals $\bm{y}-\bm{w}$ using \eqref{eq:deconv-hilbertpart}.

        {
        \begin{remark}[Choice of the Hilbert space]
            The choice of the Hilbert space $\cH$ and in particular the Hilbert norm $\norm{\cdot}_\cH$ has a significant impact on the shape of $\hfh$, involving the operator $\phih^*$ in equation~\eqref{eq:deconv-hilbertpart}. 
            Our choice of using the RKHS induced by a Gaussian kernel is motivated by the test application considered in the following Section~\ref{sec:comp:reco}, in which we want to promote more energy in the low-frequency band. 
             As demonstrated in equation~\eqref{eq:rkhs-atom}, the reconstruction elements are shifted versions of the Gaussian function $g_t$, whose spatial extension $\sigma_t$ can be tuned according to the user's needs.
        \end{remark}
        }
        
        
    \subsection{Composite reconstruction}
    \label{sec:comp:reco}

        Some algorithms directly solve Problem \eqref{eq:argmin-sps} on the continuum, for instance relying on Frank-Wolfe algorithms (Sliding Frank-Wolfe in \cite{denoyelle2019sliding}, Polyatomic Frank-Wolfe in \cite{jarret_grid_2025}, an iteratively refined Frank-Wolfe procedure in \cite{flinth2023grid}), however these methods can be long and computationally expensive to run. A fast alternative has recently been proposed in \cite{poon2025srlasso}, based on the estimation of the distance of the source locations to the knots of a uniform grid.
        To put the emphasis on the decoupling of the solutions, we decide here to simplify the solving procedure by discretizing the foreground component estimation problem~\eqref{eq:deconv-banachpart}.
        
        We do so by restricting the searched positions $z_k$ to live on a uniform fine grid. This approach alleviates the challenge of nonlinearity in the position parameters, turning \eqref{eq:deconv-banachpart} into a tractable finite dimensional generalized LASSO problem -- which can be large depending on the discretization interval and space dimensions. This strategy is usually adopted to approximate continuous-problem solutions \cite{simeoni2020functional,Debarre2019} and convergence toward gridless solutions has been studied (see \cite{duval2017sparseI,debarre2022part2} and the recent article \cite{guillemet2025a}). 

        After discretization as detailed in Appendix~\ref{app:discretization}, we obtain the following finite-dimensional problem which approximates a continuous-domain solution of \eqref{eq:deconv-banachpart}:
        \begin{equation}
            \label{eq:gridbased-sparse}
            \underset{ \bm{a} \in \mathbb{R}^{J^d}}{\arg\min} \quad ( \bm{y} - \mathbf{H}\bm{a}  )^T \mathbf{M}_{\lh}^{-1} ( \bm{y} - \mathbf{H}\bm{a}  )  + \lb \| \bm{a} \|_{1},
        \end{equation}
        The $d$-dimensional grid contains $J^d$ knots and the vector $\bm{a}\in\R^{J^d}$ stores the amplitude of Dirac impulses located on the knots.
        The matrix $\mathbf{H} \in \R^{L \times J^d}$ accounts for the application of $\phib$ to grid-based Dirac impulses.
        In the simple case $d=1$, we obtain $\mathbf{H}[\ell, j] = g(x_{\ell} - z_j)$ for $1 \leq \ell \leq L $ and $1 \leq j \leq J$.
        
        The problem \eqref{eq:gridbased-sparse} can directly be solved with a proximal gradient descent algorithm (or any atomic LASSO solver such as a Frank-Wolfe algorithm). We solve this approximate problem considering a simple illustrative example with $d=1$ in what follows.

        \subsubsection{Problem simulation and parametrization}

            We simulate a composite continuous-domain ground truth signal whose components are given by 
            $$s_1^{(0)} = \sum_{k=1}^{K_f} \beta_ k \delta(\cdot - u_k) \qquad\text{ and }\qquad s_2^{(0)} = \sum_{m=1}^{K_b} \gamma_m g_b(\cdot - v_m).
            $$
            The background component is built out of weighted replicas of the background Gaussian kernel $g_b := (1/\sqrt{2\pi\sigma_b^2}) \exp{(-x^2/(2 \sigma_b^2))}$ while the foreground involves sparse point sources. The locations $(u_k)_{k=1, \dots, K_f}$ and $(v_m)_{m=1, \dots, K_f}$ are drawn with a uniform distribution over the domain $[0, 1]$. The foreground is determined with $K_f = 2$ and $\beta_k$ drawn from a uniform distribution $\mathcal{U}([1, 10])$. The background parameters are $K_b = 100$ and $\gamma_m \sim \mathcal{U}([0.5, 1.5])$ with $\sigma_b = 0.08$.
            
            The measurement operator $\phib$ is defined with $L=100$ and $\sigma=0.02$. In addition, $\sigma_{\bm{n}}$ is set such that the signal-to-noise ratio between $\bm{y}$ and $\bm{n}$ reaches $SNR_\mathrm{dB} = 10 \log_{10}\left({\lVert \bm{y} \rVert_2^2}/{\lVert \bm{n} \rVert_2^2}\right) = 20$dB. To measure and adjust the contrast between the foreground and the background components in the signal to recover, we introduce the ratio $r_{1/2}$ of the contribution of each component in the observations, defined as
            \begin{equation}
                r_{1/2} = \frac{\lVert \phib (\fb^\dagger) \rVert_2 }{\lVert \phih (\fh^\dagger) \rVert_2 }.
            \end{equation}
    
            The observations $\bm{y}$ can be simulated with exact precision using this model. Indeed, we have
            $$
            \left[\phib(s_1^{(0)}) + \phih(s_2^{(0)})\right]_\ell = \sum_{k=1}^{K_f} \beta_k g(x_\ell - u_k) + \sum_{m=1}^{K_b} \gamma_m (g * g_b)(x_\ell - v_m),
            $$
            in which measuring the background involves the closed-form convolution
            $$
            (g\ *\ g_b)(x) = \frac{1}{\sqrt{2 \pi(\sigma^2 + \sigma_b^2)}} \exp{\left(-\frac{x^2}{2 (\sigma^2 + \sigma_b^2)}\right)}.
            $$
            
\begin{figure}[t]
                \centering
                \includegraphics[width=\linewidth, trim=0 0cm 0 .5cm, clip]{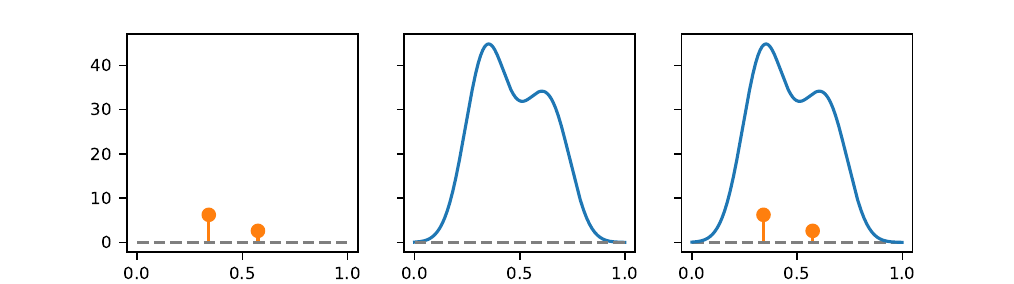}        
                \caption{Simulated source signal. \textit{Left:} Sparse component $\fb^\dagger$. \textit{Center:} Background smooth component $\fh^\dagger$. \textit{Right:} Superposition of the two components.
                }
                \label{fig:simple:source}
            \end{figure}

\begin{figure}[t]
                \centering
                \includegraphics[width=\linewidth, trim=0 0cm 0 .5cm, clip]{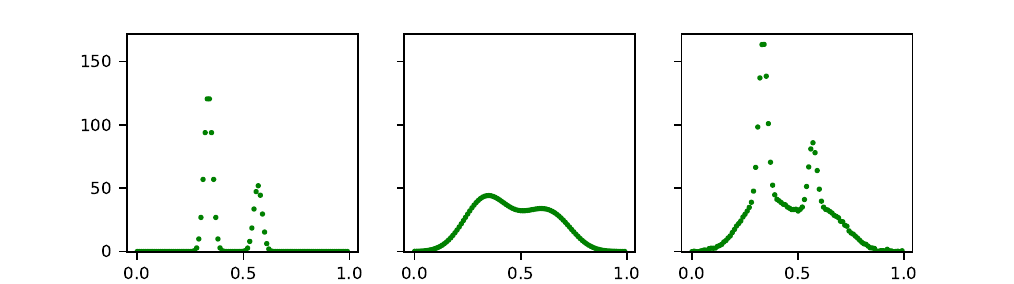}        
                \caption{Simulated measurements. \textit{Left:} Contribution of the sparse component $\phib(\fb^\dagger)\in \R^L$. \textit{Center:} Contribution of background $\phih(\fh^\dagger)\in \R^L$. \textit{Right:} Total noisy observations $\bm{y}$. In practice, only the information of the right-hand plot is accessible and respective contribution of the components is not known.}
                \label{fig:simple:measurements}
            \end{figure}

            We provide in figure~\ref{fig:simple:source} a simulated source signal built with a ratio of $r_{1/2} = 1$. 
            The associated measurements are presented in figure~\ref{fig:simple:measurements}. Note that the contributions of the two components $\phib(\fb^\dagger)$ and $\phih(\fh^\dagger)$ have the same Euclidian norm although their distribution of mass is very different, leading to the relative difference of magnitude in the right panel of figure~\ref{fig:simple:source}. 
    
        \subsubsection{Reconstruction of the signals}
        For a fine resolution reconstruction, the grid size is set to $J = n_\mathrm{srf} . L$ with the super-resolution factor fixed to $n_\mathrm{srf}=8$. The penalty parameters are tuned manually based on scaling rules to maintain values in a range consistent with the inverse problem. We set $\lh = \alpha_2 L$ with a real-valued coefficient $\alpha_2 >0$ (see Appendix~\ref{app:calculation} for a justification).
        Once $\lh$ has been set, $\lb$ is fixed as a rate of the maximum value as defined in \eqref{eq:l1max} with $\lb = \alpha_1 \lambda_{1, \mathrm{max}}$ for $0 < \alpha_1 <1$. Typically, $\alpha_1$ takes values in the range $[0.05, 0.2]$. To emphasize the smoothing effect on the background, we use the reconstruction width $\sigma_t = 0.1 > \sigma_b$. Additionally, a positivity constraint has been enforced on $s_1$ which does not break the representer theorem and usually improves the quality of the reconstruction. The approximate decoupled problem \eqref{eq:gridbased-sparse} is solved with an APGD algorithm \cite{liang2022improving}.


        \begin{figure}[t]
            \centering
            \includegraphics[width=\linewidth]{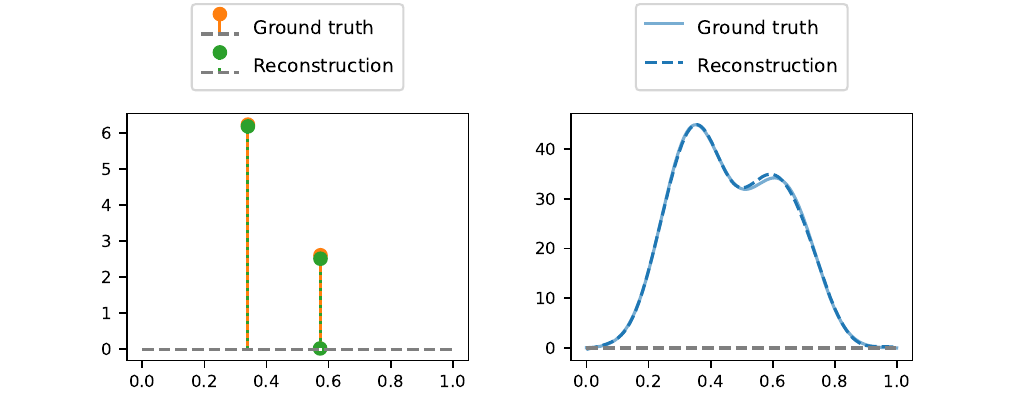}
            \caption{Recovered signals with regularization parameters $\lambda_2 = 0.4$ and $\alpha_1 = 0.05$, overlaid with the ground truth signals. \textit{Left:} Sparse foreground component. \textit{Right:} Smooth background. \textit{(Note the difference in vertical scaling.)}}
            \label{fig:simple:recos}
        \end{figure}
        

        \begin{figure}[t]
            \centering
            \includegraphics[width=\linewidth]{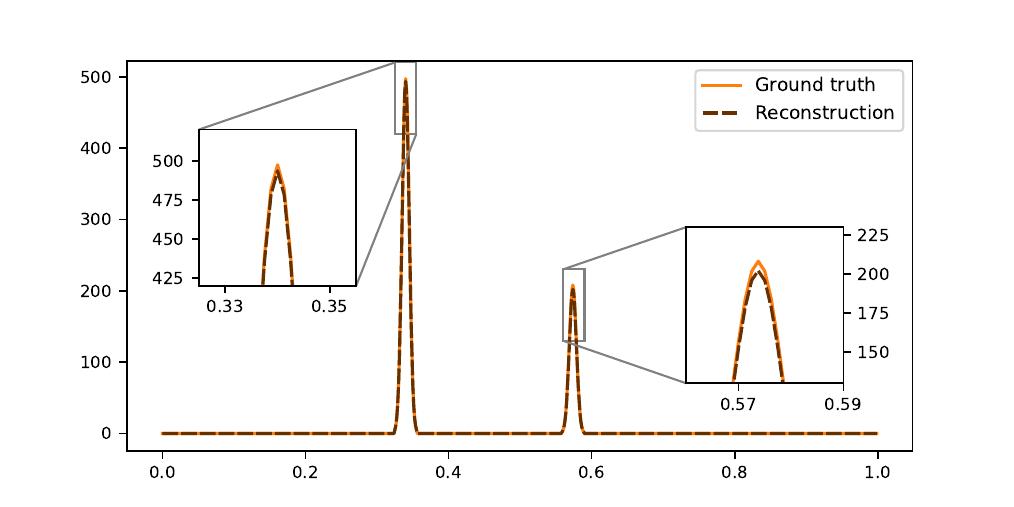}        
            \caption{Ground truth and recovered foregrounds after convolution with the sharp representation kernel, with detailed zoomed-in areas on the two impulses of the signal.}
            \label{fig:simple:recos-conv}
        \end{figure}

        Figure~\ref{fig:simple:recos} presents the foreground and background components recovered with regularization parameters $\lambda_2 = 0.04 L = 4$ and $\alpha_1 = 0.1$. In this very simple situation, we observe a strong match between the simulated source signals and the reconstructions. The recovered foreground component is composed of two clusters of impulses which accurately locate the peaks from the ground truth signal. The total intensity is however split in between the reconstruction impulses. This phenomenon is classically observed with over-parametrized LASSO problems (see for instance the discussion in \cite{duval2017sparseI}). The background reconstruction is also faithful to the ground truth signal, matching the overall maximum intensity and support of the function, even though the small variation between the two bumps is not precisely recovered.


        To provide a better comparison between the source foreground component and the recovered one, we convolve the sparse signals $\hfb$ and $s_1^\dagger$ with a representation kernel. We use a narrow Gaussian function of small standard deviation $\sigma_r = \sigma/4$. This operation intuitively blends nearby peaks while respecting the spatial spread of different clusters. The resulting signals are displayed in figure~\ref{fig:simple:recos-conv}, demonstrating consistency between the signals. Both impulses are correctly placed in the reconstruction, with a small underestimation of their intensity, which is classical side effect of $\ell_1$-type regularization.

        The regularization parameters $\lb$ and $\lh$ have a significant influence on the recovered components, and the interplay between these two parameters is not yet fully understood. To illustrate the sensitivity of the method to their respective values, Appendix~\ref{app:reconstructions} provides a grid of reconstructions obtained for various choices of regularization parameters. This analysis highlights both the stable reconstruction regimes and the parameter configurations for which the reconstruction deteriorates. Additionally, we provide in Appendix~\ref{app:more-recos} a more complex reconstruction problem, using the same simulation setup but involving $K_f = 8$ spikes in the foreground.

        \begin{remark}[About inverse crimes]
            For the background, choosing a reconstruction kernel which is different from the simulation kernel, that is, using $\sigma_t \neq \sigma_b$, prevents us from committing an inverse crime in the simulations. Although the background may be recovered with less accuracy, this situation corresponds to most application cases, in which the ground truth kernel $g_b$ is unknown.
        \end{remark}

        \subsubsection{Validation of the method}

        In addition to the reconstructed signals themselves, it is relevant to consider several indicators in order to validate the proposed composite approach, and in particular the grid-based approximation strategy employed to simplify the sparse foreground estimation problem.

\begin{figure}[t]
            \centering
            \includegraphics[width=.7\linewidth, trim=0 0.1cm 0 .6cm, clip]{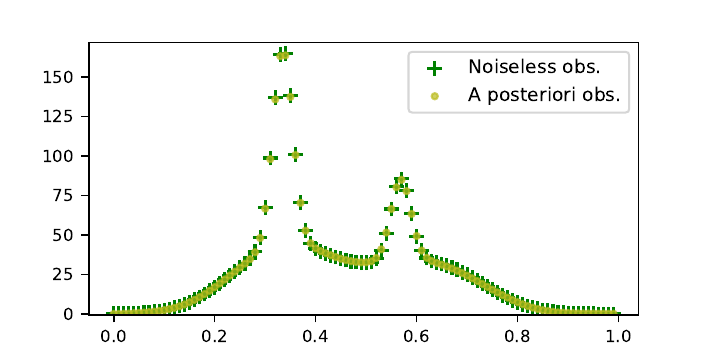}
            \caption{Comparison of the simulated noiseless observations $\{\phib(\fb^\dagger) + \phih(\fh^\dagger)\}$ (dark green crosses) and the measurements resulting from the recovered solutions $\{\phib(\hfb) + \phih(\hfh)\}$ (light green dots).}
            \label{fig:comparison_obs}
        \end{figure}

        A first classical validation consists in evaluating the mismatch between the noiseless measurements, which are known from the simulated model, and the \emph{a posteriori measurements} obtained by applying the operators $(\phib, \phih)$ to the recovered solution $(\hfb, \hfh)$. Both measurement vectors are displayed in figure~\ref{fig:comparison_obs}. We observe a strong agreement between the two datasets, with a relative error of $0.7$\% between the corresponding vectors (with respect to the Euclidian vector norm). This indicates that the recovered signals provide an accurate explanation of the simulated observations. The measurements located near the peaks of the foreground component appear slightly underestimated, which may be attributed to the classical amplitude bias encountered in B-LASSO-type problems.

        A second validation tool is obtained from the optimality conditions associated with the sparsity-promoting penalty $\cRS(\cdot) = \norm{\cdot}_\mathcal{M}$, similarly to the \emph{empirical dual certificate} commonly considered in the B-LASSO literature (see for instance \cite{duval2017sparseI,denoyelle2019sliding,jarret_grid_2025}). A pair of elements $(\fb, \fh) \in \mx \times \cHk$ is a solution of the composite problem if it satisfies the optimality conditions
        \begin{align}
            0 &\in -\phib^* (\bm{r}) + \lb \, \partial \norm{\cdot}_\mathcal{M}(\fb), \label{eq:opt:fermat}\\
            \bm{0} &= \phih^*(\bm{r}) - \lh \, \fh, \label{eq:opt:grad}
        \end{align}
        where $\bm{r} := \bm{y} - \phib(\fb) - \phih(\fh) \in \R^L$ denotes the residual vector. Equation~\eqref{eq:opt:fermat} follows from Fermat's rule using the subdifferential of the TV-norm $\partial \norm{\cdot}_\mathcal{M}(\cdot)$, while equation~\eqref{eq:opt:grad} corresponds to the vanishing of the gradient in $\cHk$ with respect to $\fh$.

\begin{figure}[t]
            \centering
            \includegraphics[width=\linewidth, trim=0 0.1cm 0 .6cm, clip]{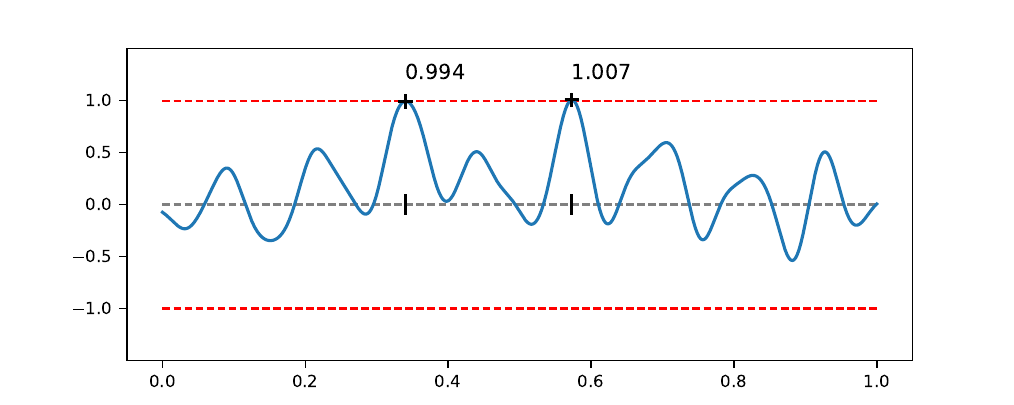}
            \caption{Empirical dual certificate $\eta(\hfb, \hfh)$ computed with the approximate solution from figure~\ref{fig:simple:recos}. The values of the local maxima are reported on the curve.}
            \label{fig:certif}
        \end{figure}

        We can then define an empirical dual certificate $\eta(\cdot, \cdot): \mx \times \cHk \to \czx$ as
        \begin{equation}
            \eta(\fb, \fh) := \frac{1}{\lb} \phib^*(\bm{y}  - \phib(\fb) - \phih(\fh)).
        \end{equation}
        For a solution pair $(\hfh, \hfb)$, this empirical dual certificate satisfies the boundedness condition $\| \eta(\hfb, \hfh) \|_\infty \leq 1$ and saturates to $\pm 1$ on the support of $\hfb$.
        In practice, the certificate is computed using the residual associated with the approximate grid-based solution obtained from problem~\eqref{eq:gridbased-sparse}. The resulting certificate is displayed in Figure~\ref{fig:certif}, together with its local maxima shown in black. The reported values fall within a $1$\% error band from the theoretical value $1$, which indicates that the grid-based solution provides a good approximation of a continuous-domain solution of problem~\eqref{eq:argmin-sps}.

        The two validation tools presented in this subsection provide quantitative indicators for assessing the sub-optimality of the proposed approach. In the present setting, they support the validity of the approximation strategy employed in the numerical implementation. More generally, they may also be used to compare different numerical configurations and to evaluate the influence of approximation parameters such as the super-resolution factor, the stopping criteria of the optimization procedures, or the implementation of the different operators.

\FloatBarrier        
\section{Benefits of our Decoupled Composite Model}
    
    Building on the simulated composite model presented in the previous section, we now illustrate the benefits of our approach.
    
    First, composite modeling accurately recovers the unknown signal in situations where single-component problems fail to distinguish foreground from background information. Second, using a decoupled numerical procedure reduces the computation time compared to a direct 2-variables approach of the composite optimization problem.

    The numerical experiments in this section have been implemented in Python based on the optimization package \texttt{Pyxu} \cite{pyxu-framework} and our code is freely accessible on the dedicated repository\footnote{\url{https://github.com/AdriaJ/decouple-composite}}. All the simulations have run on a workstation with 2 CPUs Intel Xeon E5-2680 v3 \@2.5 Ghz, 30 MB cache and 24 threads each. 

    \subsection{Compared to single-component model}
    \label{sec:bene:comp}
        When only the foreground component is of interest, we may wonder about the relevance of using a composite model in the first place. Would it be possible to recover the foreground with a single-component sparsity-promoting problem only? We address the question in this section, first visually then introducing evaluation metrics.

        \subsubsection{Single-component reconstruction}
        Using the same composite ground truth signal $(s_1^\dagger, s_2^\dagger)$ as in Section~\ref{sec:comp:reco}, we consider the following single-component B-LASSO problem with the same measurement operator $\phih$ and the same observations $\bm{y}\in\R^L$:
        \begin{equation}
            \label{eq:argmin-sparse-only}
            \underset{\fb \in \mx}{\arg\min} \frac{1}{2} \| \bm{y} - \phib (\fb) \|_2^2  + \lambda \| \fb \|_{\mathcal{M}},
        \end{equation}
        for $\lambda > 0$.
        With the fine-grid discretization proposed above, finding an approximate grid-based solution of this problem amounts to solve a classical LASSO problem. The regularization parameter $\lambda > 0$ is set specifically for this problem and is independent of $\lb$ used in \eqref{eq:deconv-banachpart}. There also exists a maximum value $\lambda_\mathrm{max} = \norm{\phib^*\bm{y}}_\infty$ for problem~\eqref{eq:argmin-sparse-only} and we set $\lambda = \alpha \lambda_\mathrm{max}$ for $0 < \alpha < 1$. It is usually larger than $\lambda_1$ as we need stronger prior information to recover a sparse solution.
        
\begin{figure}[t]
            \centering
            \includegraphics[width=\linewidth, trim=0 .8cm 0 0.4cm, clip]{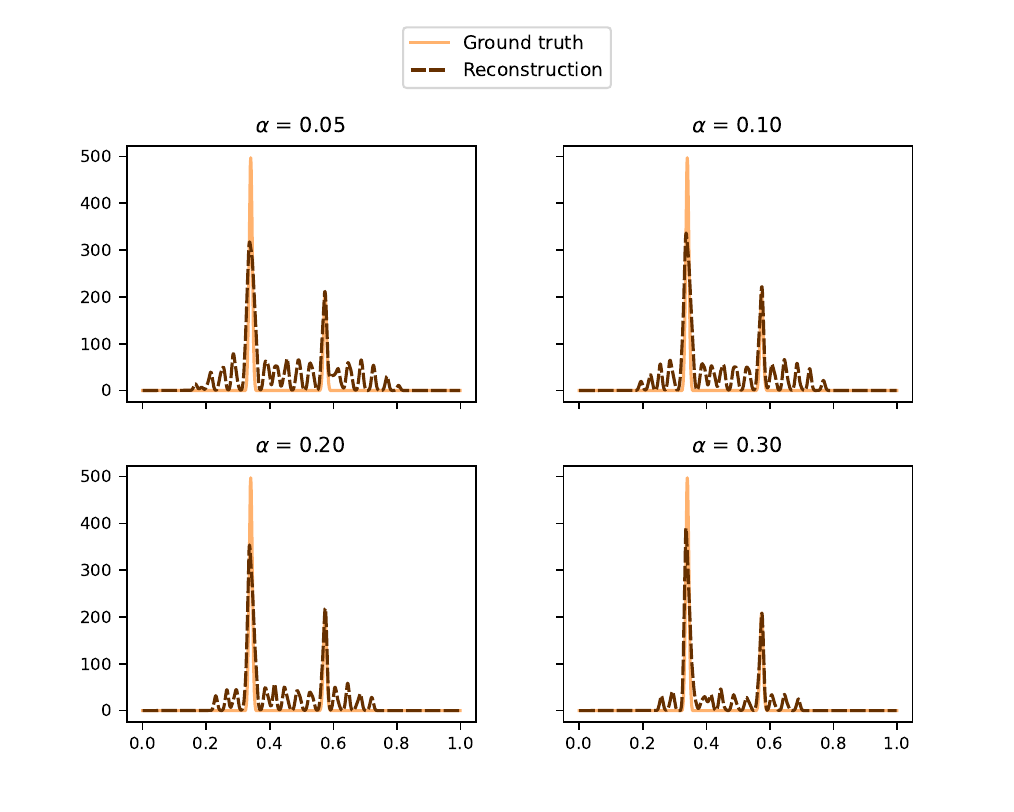}
            \caption{Single-component reconstructions using a B-LASSO problem approximated on the fine grid with several values of the regularization parameter ($\alpha=0.05, 0.1, 0.2$ or $0.3$).
            }
            \label{fig:simple:blasso-conv}
        \end{figure}

        Figure~\ref{fig:simple:blasso-conv} displays the reconstructions with the same representation kernel as before for various reconstruction parameters $\alpha$.
        We observe that the B-LASSO problem manages to recover some high intensity peaks, however the reconstructions are strongly corrupted by the presence of background information. When more sparsity is enforced with a larger value of $\alpha$, the spurious peaks tend to vanish but the intensity of the relevant peaks is also reduced, becoming lower than the source and than the composite reconstruction.


        \subsubsection{Quantitative study on the simple case}

        Assessing the quality of the reconstruction of sparse signal is always a complicated task as both intensity and localization of the recovered sources need to be simultaneously evaluated. We consider two simple metrics based on the grid-based representation of the recovered foreground. First, the ground truth and the solution are convolved with the representation Gaussian kernel. Second, the relative error is computed using either $\mathrm{L}^2$-norm or the $\mathrm{L}^1$-norm. For $\fb$ the recovered foreground component, $\fb^\dagger$ the simulated source and $g_{\sigma_r}$ the representation kernel, we define
        $$
        \mathrm{RE}_2(\fb, \fb^\dagger) = \frac{\lVert g_{\sigma_r} * (\fb - \fb^\dagger) \rVert_2 }{ \lVert g_{\sigma_r} * \fb^\dagger \rVert_2} \qquad\text{and}\qquad
        \mathrm{RE}_1(\fb, \fb^\dagger) = \frac{\lVert g_{\sigma_r} * (\fb - \fb^\dagger) \rVert_1}{\lVert g_{\sigma_r} *\fb^\dagger \rVert_1}.
        $$
        The metrics are approximated using the fine-grid representation of the signals.

        The evaluation metrics on the foreground component are reported in tables~\ref{tab:rl2} and~\ref{tab:rl1} for the sparse-plus-smooth composite model, and in table~\ref{tab:blasso} for the single-component reconstructions with the B-LASSO. A similar error table for the reconstruction of the background is provided in Appendix~\ref{app:errorbg}.

        Both metrics identify the same best pair of reconstruction parameters for the composite model, using $\alpha_1 = 0.1$ and $\lh = 8.0$. 
        The actual value of the error is significantly lower than the best error obtained with the single-component B-LASSO reconstruction. Moreover, the metrics for the B-LASSO improve with higher values of the penalty parameter, coincidentally with a decrease in intensity of the recovered signal. These observations suggest that the B-LASSO model is inappropriate for such a recovery problem and having a composite model significantly improves the reconstruction. 
        
        \begin{table}[t]
        \centering
        \begin{subtable}[t]{0.45\linewidth}
            \centering
            \begin{tabular}{l|rrrr}
                \toprule
                 & \multicolumn{4}{c}{$\alpha_1$} \\
                \cmidrule(lr){2-5}
                $\lh$ & 0.01 & 0.02 & 0.05 & 0.10 \\
                \midrule
                0.4 & 0.029 & \textbf{0.015} & 0.055 & 0.124 \\
                0.8 & 0.309 & 0.016 & 0.047 & 0.115 \\
                4.0 & 0.737 & 0.487 & 0.091 & 0.075 \\
                \bottomrule
            \end{tabular}
            \caption{Relative $\ell_2$-norm error \label{tab:rl2}}
        \end{subtable}
        \hfill
        \begin{subtable}[t]{0.45\linewidth}
            \centering
            \begin{tabular}{l|rrrr}
                \toprule
                 & \multicolumn{4}{c}{$\alpha_1$} \\
                \cmidrule(lr){2-5}
                $\lh$ & 0.01 & 0.02 & 0.05 & 0.10 \\
                \midrule
                0.4 & 0.031 & \textbf{0.015} & 0.059 & 0.134 \\
                0.8 & 0.376 & 0.018 & 0.049 & 0.124 \\
                4.0 & 1.789 & 0.975 & 0.104 & 0.077 \\
                \bottomrule
            \end{tabular}
            \caption{Relative $\ell_1$-norm error\label{tab:rl1}}
        \end{subtable}
        \caption{Evaluation metrics for the foreground component with composite reconstructions.}
        \end{table}
    
        \begin{table}[t]
            \centering
            \begin{tabular}{l|rrrr}
                \toprule
                 & \multicolumn{4}{c}{$\alpha$ in $\lambda = \alpha \lambda_{\mathrm{max}}$} \\
                \cmidrule(lr){2-5}
                 & 0.05 & 0.10 & 0.20 & 0.30 \\
                \midrule
                $\mathrm{RE}_2$ & 0.707 & 0.644 & 0.545 & \textbf{0.466} \\
                $\mathrm{RE}_1$ & 2.235 & 1.965 & 1.495 & \textbf{1.089} \\
                \bottomrule
            \end{tabular}
            \caption{B-LASSO errors \label{tab:blasso}}
        \end{table}

        \subsubsection{Error evolution with respect to the contrast}

\begin{figure}[t]
            \centering
            \includegraphics[width=\linewidth]{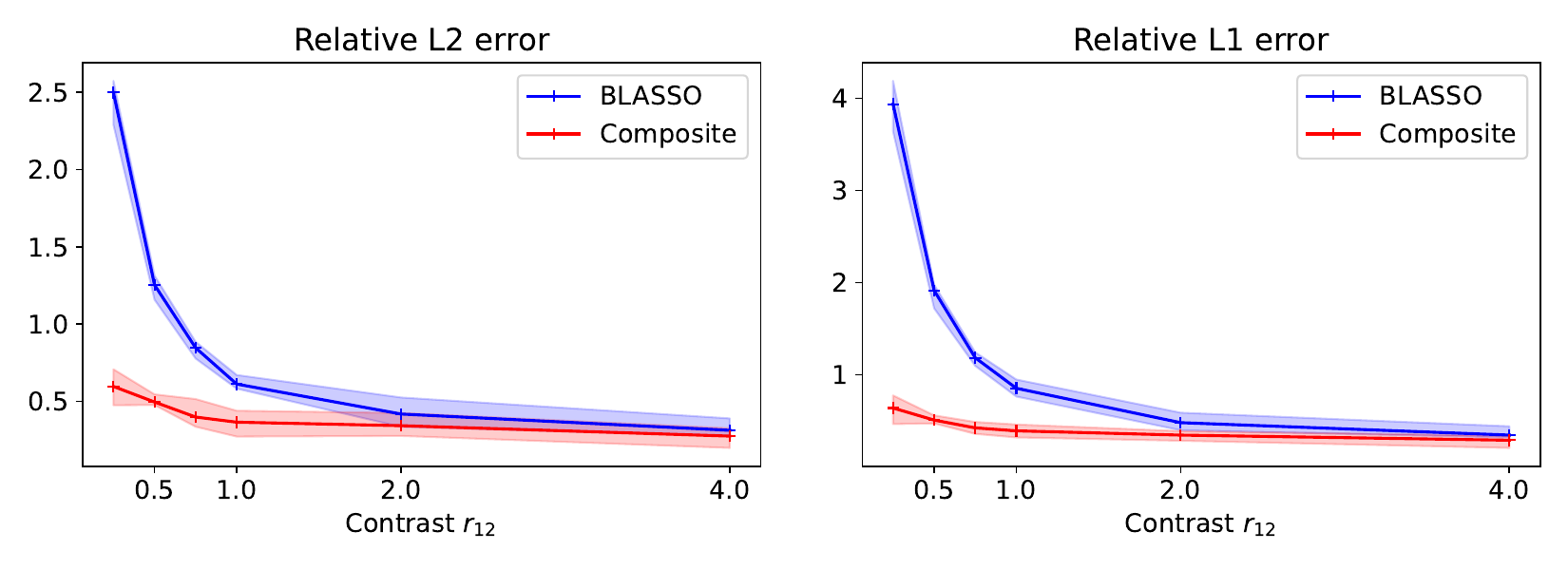}        
            \caption{Relative error with respect to contrast. For each value of $r_{1/2}$, 24 problems are simulated and solved, the median values and interquartile spreads are respectively shown with solid line and shaded area. \textit{Left:} Using $RE_2$. \textit{Right:} Using $RE_1$.}
            \label{fig:bench:metrics-vs-r}
        \end{figure}

        To further compare the interest of using a composite model over a single-component one, we study how the reconstruction metrics evolve with respect to the contrast of the signal to recover, that is the relative importance between the foreground and background components. We run a similar simulation as in Section~\ref{sec:comp:reco} while varying the parameter $r_{1/2}$ and we report in figure~\ref{fig:bench:metrics-vs-r} the best value of the metric obtained over various sets of regularization parameters. For this experiment, a more complex foreground signal has been used, with $K_f=10$, leading to overall higher error metrics.

        Independently of the contrast and the metric used, the composite model systematically outperforms the single-component one. With higher values of contrast, i.e., when the foreground is more intense relative to the background, both model reduce their error metrics and thus produce better reconstructions. Interestingly, the gap between the methods also shrinks with the contrast, ultimately being almost nonexistent for $r_{1/2} = 4$. It suggests that for sufficiently high contrasts, sparse single-component modeling may be enough to obtain an accurate reconstruction of the foreground.

    \subsection{Compared to a non-decoupled solver}
    \label{sec:bene:deco}
        So far we have illustrated Theorem~\ref{theo:main} and highlighted the interest of using a composite model in the presence of a smooth background signal. The main contribution of our theorem consists in the decoupling of the composite optimization problem and its practical benefits unveil when comparing the solving time with a direct non-decoupled approach. 

        \subsubsection{Non-decoupled approach}
        Letting apart Theorem~\ref{theo:main}, the composite optimization \eqref{eq:argmin-sps} can also be treated directly using the representer theorem of \cite{debarre2021continuous}. Indeed, it was known that there exists at least one Dirac-based sparse solution for the foreground component as
        \begin{equation}
            \hfb = \sum_{k=1}^{K_1} a_{k} \delta(\cdot - x_k)
            \label{eq:rt-deb-1}
        \end{equation}
        with $a_{k} \in \mathbb{R}$, $x_k \in \cX$ and $K_1 \leq L$.
        Additionally, the background component was known to be unique and that it can be expressed as
        \begin{equation}
            \hfh = \sum_{\ell=1}^{L} b_{\ell} \phi_\ell,
            \label{eq:rt-deb-2}
        \end{equation}
        with $b_\ell \in \mathbb{R}$ and $\phi_\ell$ the measurement functionals.

        Plugging \eqref{eq:rt-deb-1} after discretization on the fine grid $G_J$ and \eqref{eq:rt-deb-2} into the composite minimization cost function \eqref{eq:argmin-sps}, we obtain the following two-components optimization problem of dimension $J^d + L$, that we refer to as the \emph{non-decoupled approach} :
        \begin{equation}
        \label{eq:ndcp-argmin}
            \underset{(\mathbf{a}, \mathbf{b})\ \in\ \R^{J^d} \times \R^L}{\arg\min} \frac{1}{2} \| \bm{y} - \mathbf{Ha} - \mathbf{Tb} \|_2^2  + \lb \| \mathbf{a} \|_1 + \frac{\lh}{2} \langle \mathbf{b}, \mathbf{Tb} \rangle,
        \end{equation}
        where the matrix $\mathbf{T} \in \R^{L \times L}$ is defined as
        $\mathbf{T}[i, j] = (g * g_0 * g)(x_j - x_i)$ for $1 \leq i, j \leq L$
        and $\mathbf{H} \in \R^{L\times J^d}$ is from \eqref{eq:gridbased-sparse}.
        This problem is equivalent to our decoupled and discretized approach \eqref{eq:deconv-banachpart} and \eqref{eq:deconv-hilbertpart}. It is convex, finite-dimensional, and the terms are either differentiable or proximable so that the optimization can be performed with a proximal algorithm. In what follows, we solve it with APGD.
    
        \subsubsection{Quantitative assessment}

        A consequence of our representer theorem is that the composite optimization problem \eqref{eq:argmin-sps} can be solved by performing an optimization procedure on the foreground component only. We assess the practical implication of this property by comparing the runtime of a decoupled solver, including the a posteriori computation of the background component, with the non-decoupled approach of solving the two-components problem~\eqref{eq:ndcp-argmin}. For completeness, we also include in our comparison the runtime for solving a B-LASSO problem with the same input.

        To ensure a fair comparison between the solvers, the stopping criterion for all algorithms is defined as a threshold on the relative improvement between consecutive iterates of the foreground component. Moreover, the matrix–vector multiplications involving $\mathbf{H}$, $\mathbf{T}$ and $\mathbf{M}_{\lambda_2}^{-1}$ are implemented as precompiled convolutions with known kernels, so that differences in computational time accurately reflect the intrinsic complexity of the methods rather than implementation-specific optimizations.
        Similarly, both solvers rely on the APGD scheme, which is known to achieve optimal first-order convergence rates for convex problems.

\begin{figure}[p]
            \centering
            \includegraphics[width=.7\linewidth]{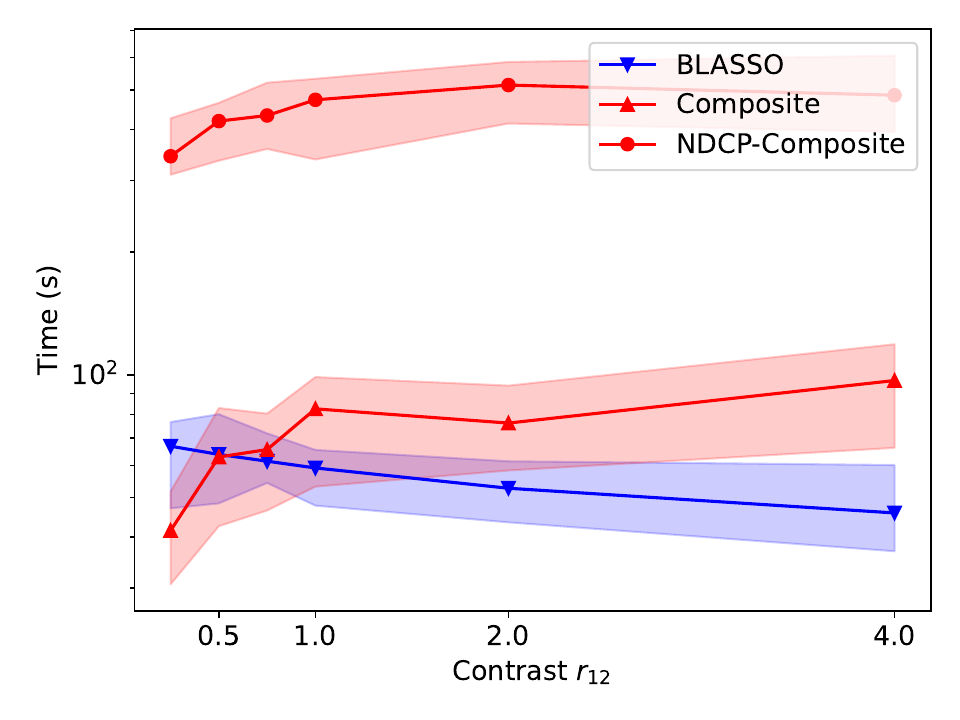}        
            \caption{Time for the best reconstruction with the different solvers with varying values of $r_{1/2}$ and fixed value of $n_\mathrm{srf} = 8$. The two composite approaches are in red, with triangle markers for the decoupled reconstruction and circle for the non-decoupled problem (``NDCP'' in the legend).}
            \label{fig:time-vs-r}
        \end{figure}

\begin{figure}[p]
            \centering
            \includegraphics[width=.7\linewidth]{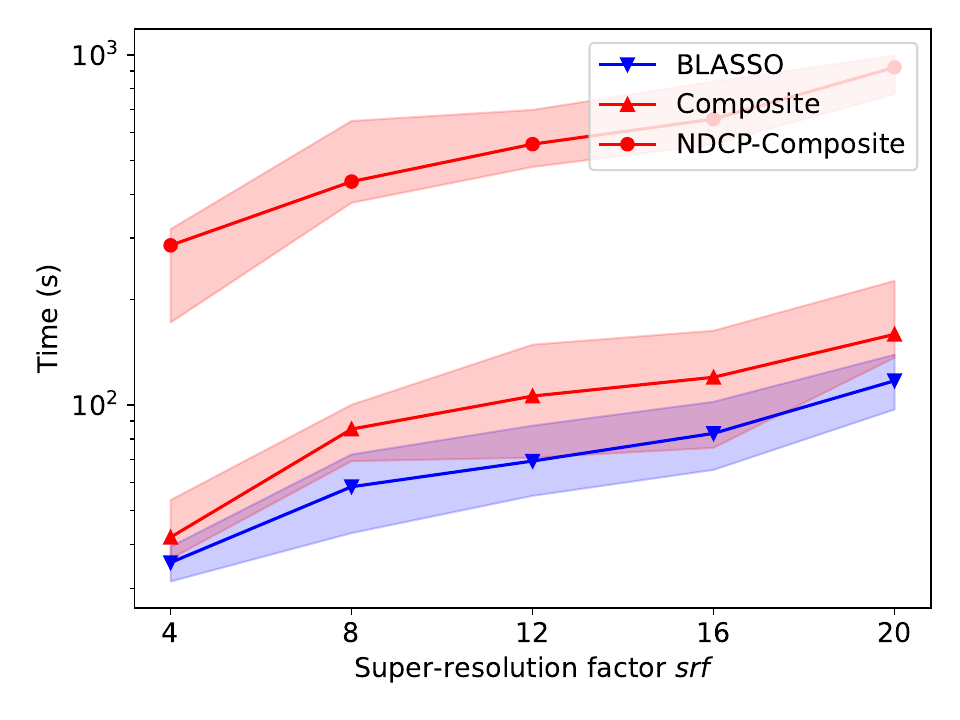}        
            \caption{Time for the best reconstruction with the different solvers with varying values of $n_\mathrm{srf}$ and fixed $r_{1/2}=1$.}
            \label{fig:time-vs-srf}
        \end{figure}

        We present two experiments: the first one reports the reconstruction time when the contrast $r_{1/2}$ evolves, in the same setup as in Section~\ref{sec:bene:comp}; the second experiment presents the reconstruction time when the super-resolution factor $n_\mathrm{srf}$ varies. The results, respectively reported in figure~\ref{fig:time-vs-r} and figure~\ref{fig:time-vs-srf}, correspond to the duration of the best reconstruction obtained through a tested set of regularization parameters. Each experiment is reproduced 24 times, the median value is reported with the solid line and the shaded area represents the interquartile spread.

        On both scenarios, the decoupled approach run significantly faster than the non-decoupled one, on average taking 15.8\% of the runtime in the first experiment (varying contrast) and 16.1\% in the second one (varying super-resolution factor). Up to numerical approximations, the solutions are identical between the two methods.
        Varying the contrast $r_{1/2}$ has little effect on the reconstruction time. Increasing the resolution, that is using more grid points in the discrete representation of $\hfb$, slows down the solving for all the algorithms. The effect is stronger for the non-decoupled method, which suffers from having its reconstruction time multiplied by approximately $3.2$.

\clearpage
\section{Conclusion}


    We introduced a new representer theorem for composite sparse-plus-smooth optimization problems, revisiting the original result from~\cite{debarre2021continuous} and partly generalizing it in a more abstract setting. Our theorem investigates deeper the connection between the two components and demonstrates a form of decoupling between them. Interestingly, the composite minimization can be transformed into an equivalent simpler single-component problem. We recover the uniqueness of the smooth component and provide a more precise closed-form expression depending on the residual of the decoupled sparse problem. Additionally, we strengthen our understanding of composite minimization problems with a general theorem on the maximum value of the regularization parameter for the Banach-penalized optimization problems.

    We highlighted the relevance of composite models for sparse recovery when the observations are corrupted with the presence of background information, in scenarios where single-component sparse modeling fails to produce accurate solutions. Moreover, the decoupled numerical procedure stemming from our representer theorem significantly outperforms the direct two-variable approach in terms of computational time.

    Building on this fast solver, composite models could be used as an enhanced version of traditional sparsity-promoting methods for practical applications involving large measurement datasets, for instance in microscopy imaging for 2D and 3D deconvolution.

\section*{Acknowledgments}
    The authors sincerely thank Martin Vetterli for his trust and guidance throughout this project. 
    A.J. is funded by the Swiss National Science Foundation (SNSF) under grant \emph{SESAM - Sensing and Sampling: Theory and Algorithms (n\textdegree 200021\_181978/1)}.
    

\clearpage
\appendix

\section{Proof of Theorem~\ref{theo:main}}
\label{app:prooftheo1}

    Let us recall the composite problem of interest. For $\lambda_1, \lambda_2 > 0$ we write the objective functional as
    \begin{equation*}
        \mathcal{J}(\fb, \fh) := \frac{1}{2} \| \bm{y} - (\phib (\fb) + \phih (\fh)) \|_2^2  + \lb \cRS(\fb) + \frac{\lh}{2} \| \fh \|_{\mathcal{H}}^2.
    \end{equation*}
    
    
    Assume first that $\fb$ is fixed and consider the optimization problem $\inf_{s_2 \in \mathcal{H}} \mathcal{J}(\fb, \fh)$. It is clearly equivalent to
    \begin{equation*}
        \underset{\fh \in \mathcal{H}}{\arg\min}  \quad \frac{1}{2} \| (\bm{y} - \phib(\fb)) - \phih(\fh) \|_2^2 + \frac{\lh}{2} \| \fh \|_{\mathcal{H}}^2,
    \end{equation*}
    whose unique solution according to Proposition~\ref{prop:hilbertRT}, depends on $\fb$ and is given by
    \begin{equation}
        \label{eq:hatfwithfb}
        \widehat{s}_{2, \fb} =\phih^* ( \phih \phih^* + \lh \mathbf{I}_L)^{-1} ( \bm{y} - \phib(\fb) ). 
    \end{equation}
    Using the latter relation, we therefore deduce that $(\hfb,\hfh) \in \mathcal{W} (\lb,\lh)$ if and only if  
    \begin{align}
        \hfh &= \widehat{s}_{2,\hfb} = \phih^* ( \phih \phih^* + \lh \mathbf{I}_L)^{-1} ( \bm{y} - \phib(\hfb) ) , \quad \text{and} \label{eq:findmethishfh}\\
        \hfb &\in\underset{\fb \in \mathcal{B}}{\arg\min}  \quad \mathcal{J}(\fb, \widehat{s}_{2, \fb}) \label{eq:newhproblem}
    \end{align}

    Replacing $\widehat{s}_{2, \fb}$ with its expression~\eqref{eq:hatfwithfb}, we observe that
    \begin{align}
        \| \bm{y} - \phih( \widehat{s}_{2, \fb} ) - \phib (\fb) \|_2^2 &= 
        \| (\mathbf{I}_L - \phih \phih^*( \phih \phih^* + \lh \mathbf{I}_L)^{-1}) ( \bm{y} - \phib (\fb) ) \|_2^2 \nonumber  \\
        &= \lh^2 \| ( \phih\phih^* + \lh \mathbf{I}_L)^{-1} (\bm{y} - \phib(\fb) ) \|_2^2 \label{eq:firstinterm}
    \end{align}
    simply using that $\mathbf{I}_L - \phih\phih^*( \phih\phih^* + \lh \mathbf{I}_L)^{-1} = \lh ( \phih\phih^* + \lh \mathbf{I}_L)^{-1}$. 

    Moreover, using again \eqref{eq:hatfwithfb} and the fact that $\phih\phih^* ( \phih\phih^* + \lh \mathbf{I}_L)^{-1} = \mathbf{I}_L - \lh ( \phih\phih^* + \lh \mathbf{I}_L)^{-1}$, the Hilbert penalty term rewrites as
    \begin{align}
        \|  \widehat{s}_{2, \fb}\|_{\mathcal{H}}^2 
        &=
        \langle \phih^* ( \phih\phih^* + \lh \mathbf{I}_L)^{-1} ( \bm{y} - \phih(\fb) ) , \phih^* ( \phih\phih^* + \lh \mathbf{I}_L)^{-1} ( \bm{y} -  \phib (\fb) ) \rangle_{\mathcal{H}} \nonumber \\
        &=
        \langle ( \phih\phih^* + \lh \mathbf{I}_L)^{-1} ( \bm{y} - \phib(\fb) ) ,  \phih \phih^* ( \phih \phih^* + \lh \mathbf{I}_L)^{-1} ( \bm{y} - \phib(\fb) ) \rangle \nonumber \\
        &= 
        \langle ( \phih \phih^* + \lh \mathbf{I}_L)^{-1} ( \bm{y} - \phib(\fb) ) ,  ( \bm{y} - \phib(\fb) ) \rangle - \lh \|( \phih \phih^* + \lh^2 \mathbf{I}_L)^{-1} ( \bm{y} - \phib(\fb) ) \|_2^2.   \label{eq:secondinterm} 
    \end{align}
    Finally, plugging the relations \eqref{eq:firstinterm} and \eqref{eq:secondinterm} together,  the cost functional $\mathcal{J}(\fb, \widehat{s}_{2, \fb})$ in \eqref{eq:newhproblem} simplifies as
    \begin{align*}
        \mathcal{J}(\fb, \widehat{s}_{2, \fb}) &=  \frac{\lh}{2} \langle ( \phih \phih^* + \lh \mathbf{I}_L)^{-1} ( \bm{y} - \phib(\fb) ) ,  ( \bm{y} - \phib(\fb) ) \rangle + \lb \cRS(\fb) \nonumber \\
        &= \frac{\lh}{2} \| ( \phih \phih^* + \lh \mathbf{I}_L)^{-\frac{1}{2}} (\bm{y} - \phib(\fb) ) \|_2^2 + \lb \cRS(\fb) \\
        &= \frac{1}{2}\| \mathbf{M}_{\lh}^{-\frac{1}{2}} (\bm{y} - \phib(\fb) ) \|_2^2 + \lb \cRS(\fb).
    \end{align*}
    This shows the first relation~\eqref{eq:banachpart}.


    Proposition~\ref{prop:abstract-rt} states that all the $\hfb$ solution of \eqref{eq:newhproblem} share the same measurement vector $\mathbf{M}_{\lh}^{-\frac{1}{2}} \phib (\hfb)$. Hence $\bm{w} = \phib(\hfb) \in \R^L$ is the common measurement vector of the Banach components $\hfb$ of the solutions $(\hfb , \hfh) \in \mathcal{W} (\lb,\lh)$, which proves \eqref{eq:hilbertpart}.

    

\section{Proof of Theorem~\ref{prop:single-lmax}}
\label{app:proof:2}

    Let us denote the objective functional $\mathcal{J}_\mathcal{B}(s) := \frac{1}{2} \norm{\bm{y} - \phib(s)}_2^2 + \lambda \norm{s}_{\cB}$. For any pair of elements $(u, s) \in \cA \times \cB$, the classical duality inequality stems from the definition of the dual norm $\norm{\cdot}_\cB$ on $\cB$ as
    \begin{equation*}
        \langle u, s \rangle_{\cA \times \cB} \leq \norm{u}_\cA \norm{s}_\cB.
    \end{equation*}
    Based on this result, we can lower-bound the value of the objective functional for any $s\in\cB$ as
    \begin{align*}
        \mathcal{J}_\mathcal{B}(s) &= \frac{1}{2} \norm{\bm{y} - \phib(s)}_2^2 + \lambda \norm{s}_{\cB} \\
                &= \frac{1}{2} \norm{\bm{y}}_2^2 + \frac{1}{2} \norm{\phib(s)}_2^2 - \langle \bm{y}, \phib(s) \rangle + \lambda \norm{s}_{\cB} \\
                &= \frac{1}{2} \norm{\bm{y}}_2^2 + \frac{1}{2} \norm{\phib(s)}_2^2 - \langle \phib^*\bm{y}, s \rangle_{\cA \times\cB} + \lambda \norm{s}_{\cB} \\
                &\geq \frac{1}{2} \norm{\bm{y}}_2^2 + \frac{1}{2} \norm{\phib(s)}_2^2 + \norm{s}_{\cB} \left(\lambda - \norm{\phib^*\bm{y}}_\cA \right) \\
                &\geq \frac{1}{2} \norm{\bm{y}}_2^2 + \frac{1}{2} \norm{\phib(s)}_2^2 + \norm{s}_{\cB} \left(\lambda - \lm \right)
    \end{align*}
    Assuming $\lambda \geq \lm$, we obtain
    \begin{equation*}
        \forall s \in \cB, \qquad \mathcal{J}_\mathcal{B}(s) \geq \mathcal{J}_\mathcal{B}(0) = \frac{1}{2} \norm{\bm{y}}_2^2,
    \end{equation*}
    hence the null element is solution, $0 \in \mathcal{V}(\lambda)$. Moreover, the strict convexity of the square $\ell_2$-norm in the data-fidelity term induces uniqueness of the fit, as detailed in Proposition~\ref{prop:abstract-rt}. For any solution $\widehat{s} \in \mathcal{V}(\lambda)$, we have 
    \begin{equation*}
        \phib(\widehat{s}) = \phib(0) = 0,
    \end{equation*}
    meaning that any solution belongs to the nullspace of the measurement operator $\phib$.
    As a consequence, they also share the same value of the penalty term (again from Proposition~\ref{prop:abstract-rt}), hence
    \begin{equation*}
        \forall \widehat{s} \in \mathcal{V}(\lambda), \qquad \norm{\widehat{s}}_\cB = \norm{0}_\cB = 0.
    \end{equation*}
    Using positive definiteness of the norm, $\widehat{s}=0$ is the unique solution in the case $\lambda \geq \lm$, which proves the first statement.

    Conversely, let us assume $0 < \lambda < \lm$ and let us prove the second statement. We want to show that $0$ is not solution, or equivalently that $\frac{1}{2} \norm{\phib(s)}_2^2 - \langle \phib^*\bm{y}, s \rangle_{\cA \times\cB} + \lambda \norm{s}_{\cB}$ can be negative. Let us write $v = \phib^*\bm{y} \in \cA$. Using the Hahn-Banach theorem \cite[Chapter~1]{Brezis_2011}, there exists an element $s_0\in\cB$ such that $\norm{s_0}_\cB = \norm{v}_\cA$ and
    \begin{equation*}
        \langle v, s_0\rangle = \norm{v}_\cA^2 = \norm{s_0}_\cB \norm{v}_\cA.
    \end{equation*}
    Using $\lm = \norm{v}_\cA$, we can then write
    \begin{align}
        - \langle v, s_0 \rangle_{\cA \times\cB} + \lambda \norm{s_0}_{\cB} &= \norm{s_0}_{\cB} \left( \lambda - \norm{v}_\cA \right) \nonumber \\
            &= \norm{s_0}_{\cB} \left( \lambda - \lm \right) < 0. \label{eq:negative}
    \end{align}
    The measurement operator $\phib: \left(\cB, \norm{\cdot}_\cB \right) \to \left(\R^L, \norm{\cdot}_2\right)$ being bounded, there exists a constant $C > 0$ such that
    \begin{equation*}
        \forall s \in \cB, \quad \norm{\phib(s)}_2^2 \leq C^2 \norm{s}_{\cB}^2.
    \end{equation*}
    Let $\xi > 0$ a real number and consider the scaled element $\xi s_0$.
    The term $\norm{\phib(\xi s_0)}_2^2$ is then of order $\xi^2$, while the positive quantity $\langle \phib^*\bm{y}, \xi s_0 \rangle_{\cA \times\cB} - \lambda \norm{\xi s_0}_{\cB}$ depends linearly on $\xi$, see equation~\eqref{eq:negative}. Consequently, as $\xi \to 0$, there exists a critical value $\xi_0 > 0$ such that, for all $0 < \xi \leq \xi_0$,  
    \begin{equation*}
        \frac{1}{2} \norm{\phib(\xi s_0)}_2^2 - \langle \phib^*\bm{y}, \xi s_0 \rangle_{\cA \times\cB} + \lambda \norm{\xi s_0}_{\cB} < 0.
    \end{equation*}
    This shows that $0$ cannot be a minimizer.

\section{Discretization of the Banach Decoupled Subproblem}
    \label{app:discretization}
        Remember that we want to solve the following problem
        $$s_1^* \in \underset{ s_1\in\mx}{\arg\min} \quad \mathcal{J}_{\mathbf{M}}(s_1)$$
        with
        $$\mathcal{J}_{\mathbf{M}}(s_1)= \|\mathbf{M}_{\lh}^{-\frac{1}{2}} ( \bm{y} - \phib(s_1)  ) \|_2^2  + \lb \| s_1 \|_{\mathcal{M}}.$$
        We know that there exist sparse solutions, which takes the form $s_1^* = \sum_k \alpha_k \delta_{z_k}$ for $z_k\in\mathcal{X}$.

        Let us introduce the $d$-dimensional fine grid $G_J$ of size $J^d$ defined as
        \begin{equation*}
            G_J = \left\{ \left(\frac{j_1}{J-1}, \dots, \frac{j_d}{J-1} \right) : 0 \leq  j_1, \dots, j_d \leq J-1 \right\}.
        \end{equation*}
        To maintain super-resolution with respect to the sampled data, the grid size $J$ is chosen much larger than the measurements grid size $K$. We define the space of Radon measures with support on this grid as
        \begin{equation*}
            V^1_J = \left\{ s \in \mathcal{M}(\mathbb{R}) : \mathrm{Supp}(s) \in G_J \right\}.
        \end{equation*}
        We then approximate the decoupled minimization \eqref{eq:deconv-banachpart} with
        \begin{equation}
            \label{eq:discrete-banach}
            s_{1, J}^* \in \underset{ s_1\in V^1_J}{\arg\min} \quad \mathcal{J}_\mathbf{M}(s_1)
        \end{equation}
        For any $s \in V^1_J$, we can write $s = \sum_j \alpha_j \delta (\cdot - z_j) $, hence
        $$
        \bm{\Phi}(s)[\ell] = \sum_j \alpha_j g(x_{\ell} - z_j)
        $$
        so that we can express 
        $$
        \bm{\Phi}(s) = \mathbf{H}\bm{a}.
        $$
        The matrix $\mathbf{H} \in \R^{L \times J^d}$ performs the convolution between the weights $\bm{a}$ and the measurement kernel $g$ sampled on the reconstruction grid and shifted to the sampling locations $x_{\ell}$. 

        Problem \eqref{eq:discrete-banach} is then equivalent to the finite-dimensional LASSO problem
        \begin{equation*}
            \underset{ \bm{a} \in \mathbb{R}^{J^d}}{\arg\min} \quad ( \bm{y} - \mathbf{H}\bm{a}  )^T \mathbf{M}_{\lh}^{-1} ( \bm{y} - \mathbf{H}\bm{a}  )  + \lb \| \bm{a} \|_{1}.
        \end{equation*}    

\clearpage

\section{Error on the Reconstruction of the Background}
    \label{app:errorbg}

    The $\mathrm{RE}_2$ error on the background is reported in the following table. The best values are obtained with a certain trade-off between the two regularization parameters.  

    \begin{table}[h]
        \centering
        \begin{tabular}{l|rrrr}
            \toprule
             & \multicolumn{4}{c}{$\alpha_1$} \\
            \cmidrule(lr){2-5}
            $\lh$ & 0.01 & 0.02 & 0.05 & 0.10 \\
            \midrule
            0.4 & 0.017 & \textbf{0.012} & 0.040 & 0.097 \\
            0.8 & 0.112 & 0.018 & 0.028 & 0.078 \\
            4.0 & 0.666 & 0.376 & 0.066 & 0.039 \\
            \bottomrule
        \end{tabular}
        \caption{Relative $\ell_2$-norm error on the background components \label{tab:rl2-bg}}
    \end{table}

\vfill

\section{Composite Reconstructions with \texorpdfstring{$K_f=2$}{Kf 2} Spikes}
    \label{app:reconstructions}

\begin{figure}[p]
    \centering
    \includegraphics[width=.91\linewidth, trim=0 .2cm 0cm .2cm,clip]{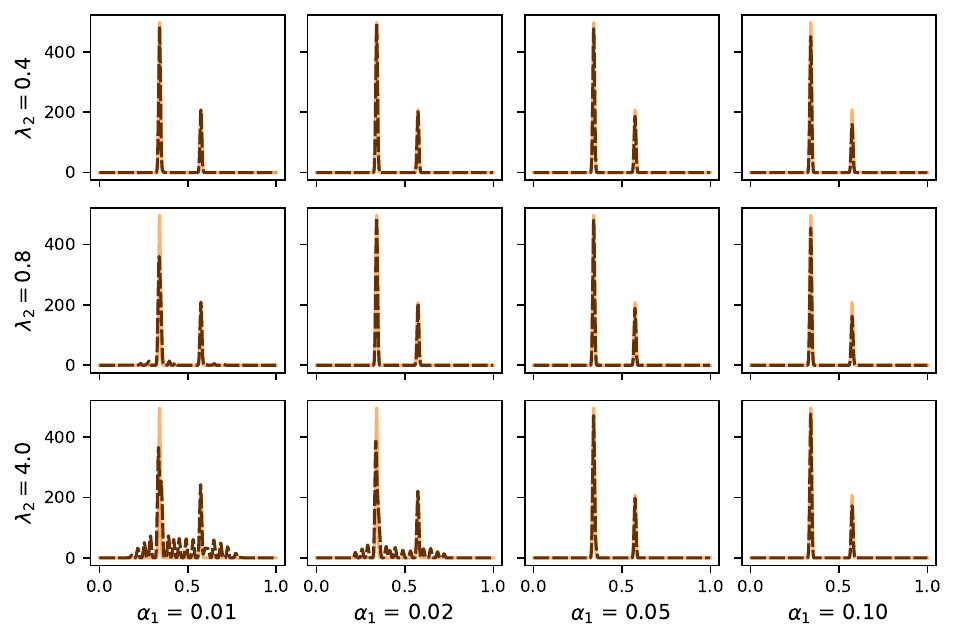}
    \caption{Foreground components after convolution with the representation kernel, for different values of
    $\lambda_2$ and $\alpha_1$. \textit{Rows :} Increasing value of $\alpha_1$ from left to right. \textit{Columns :} Increasing parameter $\lh$ from top to bottom.}
    \label{fig:simple:fg}
\end{figure}
\begin{figure}[p]
    \centering
    \includegraphics[width=.91\linewidth, trim=0 .2cm 0cm .2cm,clip]{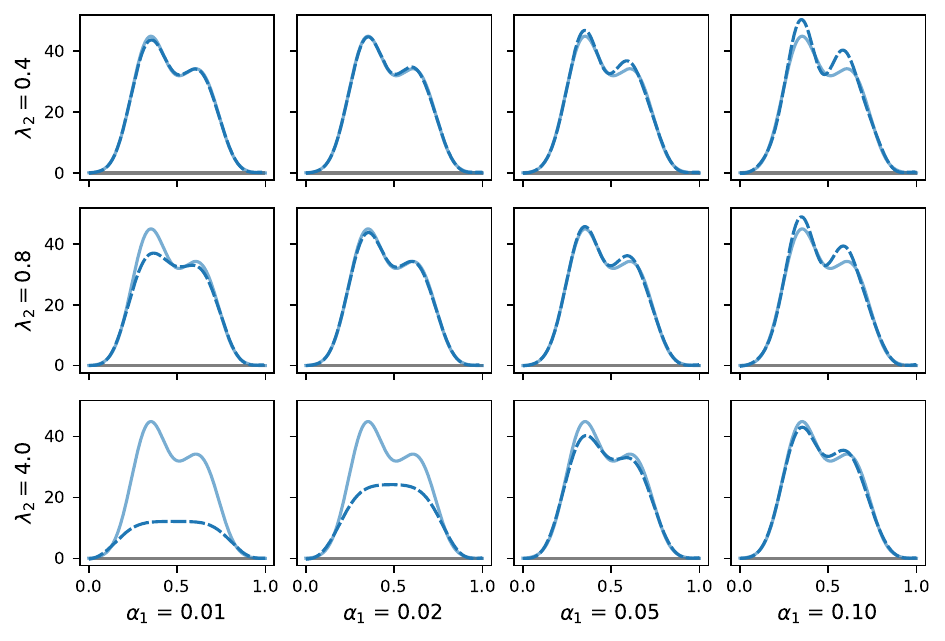}
    \caption{Background components using the same parameters.}
    \label{fig:simple:bg}
\end{figure}

    Figure~\ref{fig:simple:fg} hereafter displays the various reconstructions of the foreground component when varying the regularization parameters. Figure~\ref{fig:simple:bg} does the same for the background component. Note that the value of $\lb$ is not provided in these plots, only the coefficient $\alpha_1 = \lambda_1 / \lambda_{1, \mathrm{max}}$. This reparametrization somehow interleaves the value of the parameters, as $\lambda_1$ strongly depends on the value of $\lambda_2$. Consequently, in each column, $\lambda_1$ varies across the rows. 

    \vfill

\clearpage

\vfill

\section{Composite Reconstructions with \texorpdfstring{$K_f=8$}{Kf 8} Spikes}
    \label{app:more-recos}

    \vfill

\begin{figure}[p]
        \centering
        \includegraphics[width=.91\linewidth, trim=0 .2cm 0cm .2cm,clip]{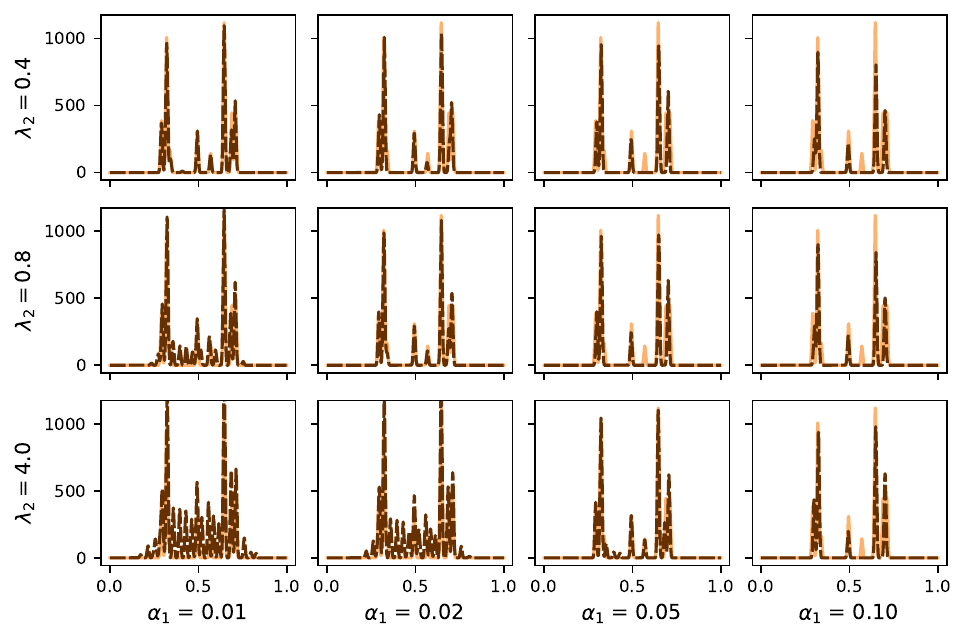}
        \caption{Foreground components after convolution with the representation kernel, for different values of
        $\lambda_2$ and $\alpha_1$. \textit{Rows :} Increasing value of $\alpha_1$ from left to right. \textit{Columns :} Increasing parameter $\lh$ from top to bottom.}
        \label{fig:complex:fg}
    \end{figure}

\begin{figure}[p]
        \centering
        \includegraphics[width=.91\linewidth, trim=0 .2cm 0cm .2cm,clip]{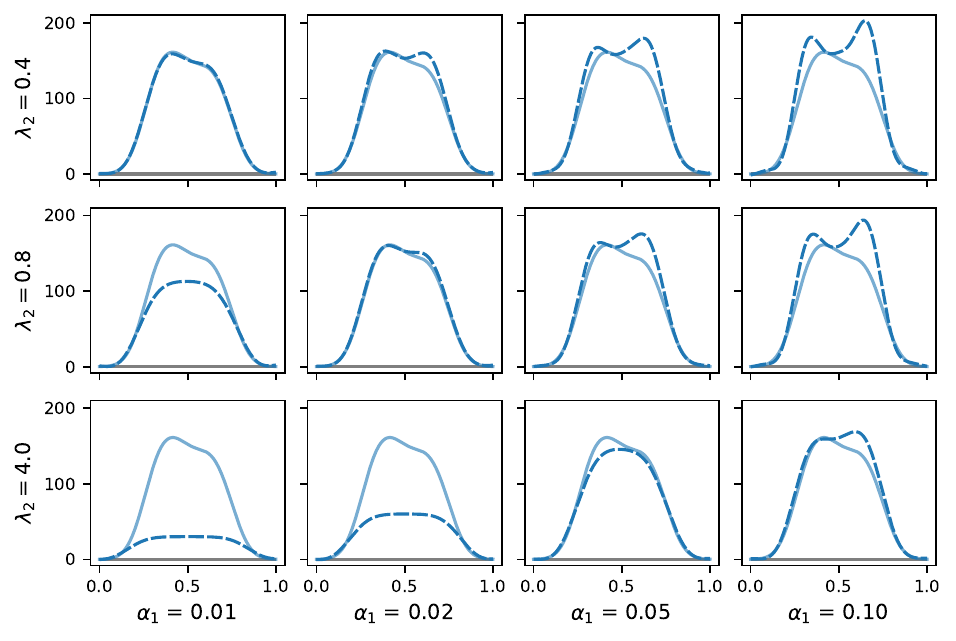}
        \caption{Background components using the same parameters.}
        \label{fig:complex:bg}
    \end{figure}

    To provide more insights into the mechanics of composite reconstruction, we illustrate here the recovery performance using a more complex foreground involving $K_f=8$ spikes. The rest of the simulation setup is identical as described in Section~\ref{sec:comp:reco}. The recovered foreground and background signals are respectively displayed in figures~\ref{fig:complex:fg} and ~\ref{fig:complex:bg}. The associated metrics are presented in tables~\ref{tab:complex:rl2} and~\ref{tab:complex:rl1} for the foreground and table~\ref{tab:complex:rl2} for the background.

    We observe that overall the metrics deteriorate compared to the simpler example with $K_f=2$, more particularly for the foreground component. The presence of more impulses increases the potential locations of mismatch in the recovered intensity. However, the best reconstruction, obtained with $\alpha_1 = 0.05$ and $\lh=0.4$, still manages to recover a visually satisfying signal. Even closely located peaks or smaller intensity ones can be distinguished. In this low-noise regime ($SNR_\mathrm{dB}=20$dB), composite reconstruction is able to decompose the foreground from the background, which is promising for real-world applications.

    \vfill

    \begin{table}[h]
    \centering
    \begin{subtable}[t]{0.45\linewidth}
        \centering
        \begin{tabular}{l|rrrr}
            \toprule
             & \multicolumn{4}{c}{$\alpha_1$} \\
            \cmidrule(lr){2-5}
            $\lh$ & 0.01 & 0.02 & 0.05 & 0.10 \\
            \midrule
            0.4 & \textbf{0.185 }& 0.425 & 0.539 & 0.528 \\
            0.8 & 0.430 & 0.241 & 0.551 & 0.535 \\
            4.0 & 0.745 & 0.591 & 0.333 & 0.582 \\
            \bottomrule
        \end{tabular}
        \caption{Relative $\ell_2$-norm error \label{tab:complex:rl2}}
    \end{subtable}
    \hfill
    \begin{subtable}[t]{0.45\linewidth}
        \centering
        \begin{tabular}{l|rrrr}
            \toprule
             & \multicolumn{4}{c}{$\alpha_1$} \\
            \cmidrule(lr){2-5}
            $\lh$ & 0.01 & 0.02 & 0.05 & 0.10 \\
            \midrule
            0.4 & \textbf{0.175} & 0.472 & 0.609 & 0.574 \\
            0.8 & 0.654 & 0.250 & 0.622 & 0.588 \\
            4.0 & 1.299 & 1.000 & 0.326 & 0.660 \\
            \bottomrule
        \end{tabular}
        \caption{Relative $\ell_1$-norm error\label{tab:complex:rl1}}
    \end{subtable}
    \caption{Evaluation metrics for the foreground component presented in figure~\ref{fig:complex:fg}.}
    \end{table}

    \vfill

    \begin{table}[h]
        \centering
        \begin{tabular}{l|rrrr}
            \toprule
             & \multicolumn{4}{c}{$\alpha_1$} \\
            \cmidrule(lr){2-5}
            $\lh$ & 0.01 & 0.02 & 0.05 & 0.10 \\
            \midrule
            0.4 & \textbf{0.014} & 0.083 & 0.179 & 0.286 \\
            0.8 & 0.260 & 0.034 & 0.152 & 0.248 \\
            4.0 & 0.782 & 0.581 & 0.085 & 0.118 \\
            \bottomrule
        \end{tabular}
        \caption{Relative $\ell_2$-norm error on the background components presented in figure~\ref{fig:complex:bg} \label{tab:complex:rl2-bg}}
    \end{table}

    \vfill

\clearpage
\section{Calculation for the Experiments}
    \label{app:calculation}

    \subsection{Computation of \texorpdfstring{$M_{\lh}$}{M lambda 2}}

    Using the definition of the PSF $g$ in Section~\ref{sec:application}, we derive the expression of the matrix $\mathbf{M}_{\lambda_2}$ defined in equation~\eqref{eq:def-Mphi}. For $1\leq k, \ell \leq L$, we have:
    \begin{equation*}
        \mathbf{M}_{\lh}[k, \ell] = \frac{1}{\lh} \left( \langle \phihh_k , \phihh_\ell \rangle_{\mathcal{H}} + \lh \delta[k - \ell] \right)
    \end{equation*}
    Let us compute the inner product between the functionals $\phihh$
    \begin{equation*}
        \begin{aligned}
            \langle \phihh_k , \phihh_\ell \rangle_{\mathcal{H}} &= (g * \phihh_\ell)(x_k)\\
                &= \left(g * \left[(g_0 * g)(x_\ell - \cdot)\right]\right)(x_k) \\
                &=  (g * g_0 * g)(x_\ell - x_k)
        \end{aligned}
    \end{equation*}
    using equation~\eqref{eq:defPhihStar} and the symmetry of the kernels $g$ and $g_0$.
    We finally obtain the expression
    \begin{equation*}
        \label{eq:mlambda-1d}
        \mathbf{M}_{\lambda_2}[k, \ell] = \frac{1}{\lh} \left( (g * g_0 * g)(x_\ell - x_k) + \lh \delta [k - \ell] \right) .
    \end{equation*}

    \subsection{Scaling of \texorpdfstring{$\lh$}{lambda 2}}
    From equation~\eqref{eq:def-Mphi}, the matrix $\mathbf{M}_{\lambda_2}$ can be written as
    \begin{equation}
        \label{eq:mlambda2}
        \mathbf{M}_{\lh} = \mathbf{I}_L +  \frac{1}{\lh}\phih\phih^*
    \end{equation}
    For any vector $\mathbf{h} \in \R^L$, we have
    \begin{align*}
        \forall 1 \leq k \leq L, \quad \left(\phih\phih^*\mathbf{h}\right)[k] &= \sum_\ell \langle \phihh_k , \phihh_\ell \rangle_{\mathcal{H}} h_\ell \\
            &= \sum_\ell u_{k - \ell} h_\ell,
    \end{align*}
    with $u_i = (g * g * g_0 * g_0)(x_i)$ for $ \leq i \leq L$. Assuming the support of $\mathbf{h}$ is concentrated in the center of the vector and there is no information near the borders, $\phih\phih^*\mathbf{h}$ can be seen as a convolution between $\mathbf{h}$ and $\mathbf{u}$.

    It is possible to prove that the Lipschitz constant of the operator $\phih\phih^* \in \R^{L \times L}$, which is its maximum singular value, is bounded by $\sqrt{L} \norm{\mathbf{u}}$. Using that $\norm{g*g*g_0*g_0}_2 = 1$, we can approximate $\norm{\mathbf{u}}_2 \approx \sqrt{L}$ using Riemann sum. Hence, we approximate the maximum singular value of $\phih\phih^*$ with $L$.

    To maintain the term $\frac{1}{\lh} \phih\phih^*$ of the same order as $\mathbf{I}_L$ in $\eqref{eq:mlambda2}$, we set $\lh$ as a ratio of $L$ and we obtain our proposed scaling rule $\lh = \alpha_2 L$.
    
\clearpage

\bibliographystyle{vancouver}
\bibliography{refs}

\end{document}